\documentclass{article}
\usepackage[utf8]{inputenc}
\usepackage{amsmath,amssymb}
\usepackage{dsfont}
\usepackage[numbers]{natbib}
\usepackage{color}

\newtheorem{theorem}{\bf Theorem}[section]
\newtheorem{proposition}{\bf Proposition}[section]
\newtheorem{lemma}{\bf Lemma}[section]
\newtheorem{remark}{\bf Remark}[section]
\newtheorem{example}{\bf Example}[section]

\def\eps{\varepsilon}
\def\R{\mathbb{R}}
\def\T{\mathbb{T}}

\renewcommand{\div}{\mathrm{div}\,}
\newcommand\transp[1]{#1^\dagger}
\newcommand{\tr}{\mathrm{tr}\,}
\newcommand{\bsigma}{\boldsymbol{\sigma}}
\newcommand{\bSigma}{\boldsymbol{\Sigma}}

\newcommand{\bdelta}{\boldsymbol{\delta}}
\newcommand{\bA}{\boldsymbol{A}}
\newcommand{\bB}{\boldsymbol{B}}
\newcommand{\bC}{\boldsymbol{C}}
\newcommand{\bQ}{\boldsymbol{Q}}
\newcommand{\bzero}{\boldsymbol{0}}
\newcommand{\finpreuve}{\hfill $\square$\\}

\title{Global strong solutions for some differential viscoelastic models}

\author{Laurent Chupin\thanks{Universit\'e Clermont Auvergne, Laboratoire de Math\'ematiques Blaise Pascal CNRS-UMR 6620, Campus des C\'ezeaux, F-63177 Aubi\`ere cedex, France ({\tt Laurent.Chupin@uca.fr})}}

\date{Version: \today}

\begin{document}

\maketitle

\begin{abstract}
The purpose of this article is to show that there are many differential viscoelastic models for which the global existence of a regular solution is possible.  Although the problem of global existence in the classic Oldroyd model is still open, we show that by adding a non-linear contribution (proposed by R.G. Larson in~1984), it is possible to obtain more regular and global solutions, regardless of the size of the data (in the two-dimensional and periodic case).
Similarly, more complex appearance models such as those related to "pom-pom" polymers are interesting and mathematically richer: some "natural" bounds on the stress make it possible to obtain global results.  On the other hand, in the last part, we show that other models clearly do not seem to fit into this framework, and do not even seem to have a global solution in time.
These kinds of results allow to highlight the advantages and disadvantages of such or such viscoelastic fluid models. They can thus help rheologists and numericists to make choices with new arguments.
\end{abstract}





\section{Introduction}\label{introduction}

Throughout this paper we will consider incompressible fluids.
We model the flow using partial differential equations derived from fluid mechanics principles. More specifically, denoting by~$v$ the fluid velocity field,~$\rho$~its density and by~$\bSigma$ the Cauchy stress tensor, the fundamental law of Dynamics (Newton’s law) yields the equation
\begin{equation}\label{NS}
\left\{
\begin{aligned}
& \rho(\partial_t v + v\cdot \nabla v) = \div \bSigma + f\\
& \div v = 0,
\end{aligned}
\right.
\end{equation}
where~$f$ denotes some external forces applied to the fluid.
When you consider a Newtonian solvent in dilute polymer solution, it is usual to write the Cauchy stress tensor as the sum of three main contributions corresponding to the pressure forces, the viscous effects and the elastic ones:
\begin{equation}\label{NS1}
\bSigma = -p\, \bdelta + \eta_s Dv + \bsigma.
\end{equation}
The scalar $p$ is the hydrostatic pressure,~$Dv$ is the rate of deformation tensor, $Dv = \frac{1}{2} (\nabla v + \transp{(\nabla v)})$, $\eta_s$ is the kinematic viscosity of the solvent, and $\bsigma$ corresponds to the extra-stress tensor.
While~\eqref{NS} is derived from first principles and is thus general, the precise form of the extra-stress tensor~$\bsigma$ is given by a constitutive law which depends on the fluid behavior.
For some fluids like polymer melts, biological fluids, etc. the value of the stress tensor at the present time~$t$ depends on the history of the past deformations (of course with a weaker dependence on the far past), and not only on the present deformation. In other words, such fluids have memory and property of elasticity.

\paragraph{A wide variety of models.}

The model which is at the base of the models presented in this article is the nonlinear one described by Oldroyd, see for instance~\cite{Oldroyd50}. The Oldroyd-B fluid presents one of the simplest constitutive models capable of describing the viscoelastic behavior of dilute polymeric solutions under general flow conditions. He postulated quasi-linear and nonlinear
constitutive equations of differential and integral types related to the external observable variables, the extra-stress tensor~$\bsigma$, and strain rate tensor~$Dv$, and also elucidated some of the important principles of invariance.
The most iconic model of Oldroyd's work remains the following relatively simple model (called the Oldroyd-B model):
\begin{equation}\label{eq:Oldroyd}
\lambda \big( \partial_t \bsigma + v\cdot\nabla \bsigma - \bsigma\cdot \nabla v - \transp{(\nabla v)}\cdot \bsigma \big) + \bsigma  = 2 \eta_p Dv
\end{equation}
where $\lambda$ is the relaxation time and~$\eta_p$ corresponds to the polymer viscosity.\par
Several other concepts were also developed by scientists such as Rivlin, Green, Tobolsky, Ericksen, Lodge, Phan-Thien, Tanner, Giesekus, Doi, Edwards and their numerous successors (see the references~\cite{Bird77,Doi88,Lodge64,Porod66,Rivlin65,Saut13,Tanner77}). It was also recognized that many of these concepts were associated with the Oldroyd approach.
The results that are discussed in this article only relate to such differential models.

\begin{remark}
Note that there are other classes of models - that are not discussed in this article, for instance:
\begin{itemize}
\item Models of integral type for which the extra-stress tensor in one time $t$ is expressed using an integral covering all the past times $T<t$ of the stresses/strain.
The best-known models of this type are included in the K-BKZ model class, which the review of Mitsoulis~\cite{Mitsoulis13} gives an excellent overview. See also~\cite{KBKZ2} or~\cite{Saut13}.
\item Micro-macros models, coupling the mesoscopic scale of kinetic theory to the macroscopic scale of continuum mechanics. Generally, the extra-stress takes the form of a particular average computed with the distribution function over all possible configurations, the distribution function associated being the solution to a Fokker-Planck equation, see for instance~\cite{LeBris12,Ottinger12,Saut13}.
\end{itemize}
\end{remark}

\paragraph{But few mathematical global results.}

About the Oldroyd-B model, the question of existence of global solution is still open, even in the two dimensional case. There are, nevertheless, partial results: Guillopé and Saut~\cite{Guillope-Saut4,Guillope-Saut-CRAS,Guillope-Saut3,Guillope-Saut1} proved the existence of local strong solutions; Fernández-Cara, Guillén and Ortega~\cite{Fernandez-Guillen-Ortega-CRAS,Fernandez-Guillen-Ortega} proved local well posedness in Sobolev spaces; Chemin and Masmoudi~\cite{Chemin} proved local well posedness again but in critical Besov spaces.
In these papers, some global existence results hold, only assuming small data. It may be noted that in the co-rotational case, that is using a peculiar time derivative for the stress tensor, see the definition~\eqref{derive1} hereafter, and only for that case, a global existence result of weak solution has been shown, see the result of Lions and Masmoudi~\cite{Lions-Masmoudi-viscoelastique}.\\
Concerning the micro-macro models, coupling the Navier-Stokes equations with a Fokker-Planck equation, a lot of local existence results are proved, see for instance~\cite{Jourdain-Lebris-Lelievre,Masmoudi08,Renardy,Zhang-Zhang}. Recently, Masmoudi~\cite{Masmoudi13} proved global existence of weak solutions to the FENE (Finite Extensible Nonlinear Elastic) dumbbell model.\\
The integral models are also studied in the last fifty years.
The first significant theoretical results are probably due to Kim~\cite{Kim}, M. Renardy~\cite{Renardy3}, Hrusa and Renardy~\cite{Renardy2}, Hrusa, Nohel and Renardy~\cite[Section IV.5]{Renardy1}.
In all these works the solution is either local in time or global but with small data. Later, Brandon and Hrusa~\cite{Brandon-Hrusa} study a one dimensional model with a singularity in the non-linearity: they obtain global existence results for sufficiently small data.
Recently global existence results for integral law are given in~\cite{Chupin14}.\\
One of the objectives of the present article is to show that there exists other realistic viscoelastic fluid models for which global solution existence results can be shown. This will provide additional arguments for choosing models.

\paragraph{The PEC model as an example.}

In the present article, we will first focus on a typical model but a discussion of the relevance of identical results for other more or less complex models is proposed in the last sections. All models will be written in dimensionless form and we have chosen to take all the constant equal to~$1$ (Reynolds number, Weissenberg number, retardation parameter...).
In practice, all the results announced will be true with any coefficients, the proofs remaining similar.\\

The model studied was initially introduced by Larson in~\cite{Larson84} to take into account the partially extensible nature of the polymer strands (without extensibility, we get the Maxwell model whereas the complete extensibility model corresponds to the Doi-Edwards model). This model, called PEC model, writes\footnote{See the Appendix on page~\pageref{sec:notations} for the notations, especially for contractions $\cdot$ or $:$ and for transposition $\transp{}$.}
\begin{equation}\label{PEC1}
\left\{
\begin{aligned}
& \partial_t v + v\cdot\nabla v - \Delta v + \nabla p = \div \bsigma \\
& \div v = 0 \\
& \partial_t \bsigma + v\cdot\nabla \bsigma - \bsigma\cdot \nabla v - \transp{(\nabla v)}\cdot \bsigma + \eps \, (Dv : \bsigma) \bsigma = \bzero.
\end{aligned}
\right.
\end{equation}
The trilinear term $\eps\, (Dv : \bsigma) \bsigma$ makes this model a singular case and truly different from more classical models like the Oldroyd model, see the equation~\eqref{eq:Oldroyd}. In a relatively surprising way, we will see that this term has a fundamental role for the stress as soon as $\varepsilon$ satisfies $0<\varepsilon<\frac{2}{\tr \bsigma|_{t=0}}$.\\

For this model, we arrive at a result of the existence, uniqueness and stability of a global strong solution, without smallness assumption on the data.
This result, as all those presented in this paper, is proved in the two-dimensional case, and assuming that the domain is periodic in space. The dimensional condition is constrained by the known results on the Navier-Stokes equations: it is well known that currently we do not know how to show the existence of a strong global solution to the Navier-Stokes equations in the three-dimensional case: it is one of the Millennium Problems. In contrast, the periodicity condition is a technical condition that allows overcoming some potential boundary problems: It is likely that the results can be improved by avoiding this condition. Right now, we don't see how we can do that.
More precisely, we show (see Theorem~\ref{th:global-existence-1} on page~\pageref{th:global-existence-1}):

\begin{theorem}\label{th:intro}
Given any regular initial data for the velocity~$v|_{t=0}$ and for the stress~$\bsigma|_{t=0}$, the system~\eqref{PEC1} admits a unique strong solution $(v,p,\bsigma)$ defined for all time $t>0$ (in the two-dimensional case). Moreover, solution's behavior changes continuously with the initial conditions.\par
In the three dimensional case, the result holds if we do not taking into account the convective term $v\cdot \nabla v$ in the velocity evolution equation.
\end{theorem}

Note that, as usually, the uniqueness of the pressure is obtained by adding a condition, for example of zero average pressure.

\paragraph{Why can we prove a global existence result for these kinds of models?}

The key point of the proofs is the obtention of bounds for the stress (see Lemma~\ref{lem:stress}). To obtain these bounds, we will see that this model can be written using a conformation tensor as unknown instead of the constraint tensor (see Subsection~\ref{subsec:PEC-formulation}). The bound is then relatively natural.
Once we have precise information on stress, we can then use fine properties on the Navier-Stokes equations, see Lemma~\ref{lem:velocity}, to deduce limits on the velocity field. In this lemma, the estimate~\eqref{PEC33} plays a pivotal role in the argument. It is this point which is the main difference to the one-dimensional case studied in~\cite{Renardy1}, where a much simpler estimate holds.
The last step consists in showing that these bounds prohibit the finite time explosion of solutions (mainly via a Gronwall estimate).

\paragraph{Notion of "good" solution in the mathematical sense}

In 1902, J. Hadamard~\cite{Hadamard} introduced the notion of well-posed problem. He proposed that mathematical models of physical phenomena should have the following three properties:\\
\phantom{toto} 1/ a solution exists;\\
\phantom{toto} 2/ it is unique;\\
\phantom{toto} 3/ its behavior changes continuously with the initial conditions.\\
Clearly, Theorem~\ref{th:intro} previously announced implies that the PEC problem~\eqref{PEC1} is well posed in the Hadamard sense.
Moreover, for an evolution model, i.e. a time-dependent model, it is possible to add another condition:\\
\phantom{toto} 4/ the solution exists for all future times.\\
The result announced by the Theorem~\ref{th:intro} indicates that the PEC model described by~\eqref{PEC1} has this property, and one of the objectives of this article is not only to prove this result but also to show that other viscoelastic models have this property.\\
It should be noted that, generally speaking, there are models for which the existence of a solution is only local in time. This fact may in some cases correspond to a physical reality, for example when the model ceases to be valid.
However, a model can be described as "bad" if one of its physical properties is no longer satisfied. This will be the case in the example described in Subsection~\ref{sec:false} where we will exhibit viscoelastic fluid models for which the conformation tensor does not remain positive.

\paragraph{Article outline}
The paper is organized as follows:
The main section of the paper is the next section (Section~\ref{sec:PEC}) in which we precisely present the PEC model, we explicitly introduce the main results and we give all the mathematical proofs.
In Section~\ref{sec:others} we discuss the relevance and possible adaptability of the proof to several other differential models for polymer flows: PTT, MGI, Pom-Pom...
Section~\ref{sec:conclusion} is a conclusion. It also contains example of model for which the previous study would not be successful. An additional section - Appendix - indicates the main notations used in the paper.

\section{The partially extending strand convection (PEC) model}\label{sec:PEC}

\subsection{Structure and origin of the PEC model}

The PEC model is based on molecular theory,
for a polymeric melt composed of polymers with side branches, which hinder full retraction.
This model was originally introduced by Larson~\cite{Larson84}. The idea behind this differential model is that long side branches can limit the extent to which a polymer strand is able to contract back within the confining tubelike region postulated by Doi and Edwards.
More precisely, on page~554 of~\cite[equation (41)]{Larson84}, the relation between~$\bsigma$ and~$Dv$ reads
\begin{equation}\label{PEC3}
\overset{\bigtriangledown}{\bsigma} + \frac{2\xi'}{3G} (Dv : \bsigma) \bsigma = \bzero,
\end{equation}
where $\xi'\in (0,1)$ and $G>0$ are two physical constants.
As indicated in~\cite{Larson84}, constant~$\xi'$ can be seen as an interpolation parameter between the Lodge model (corresponding to the case $\xi'=0$) and the Doi-Edwards model (corresponding to the case $\xi'=1$).
We also note that in order to take into account the reptational diffusion, as well as extending convection, it is possible to introduce the following constitutive equation (see~\cite[equation (54), page 557]{Larson84}):
\begin{equation*}\label{PEC4}
\tau_0 \Big( \overset{\bigtriangledown}{\bsigma} + \frac{2\xi'}{3G} (Dv : \bsigma) \bsigma \Big) + \bsigma = G \bdelta,
\end{equation*}
where the linear term arises from the Doi-Edwards' or impulsive diffusion (the parameter~$\tau_0$ corresponding to a relaxation time).
To treat a more general case, we are interested in the constitutive equation where the relation between stress $\bsigma$ and velocity gradient~$Dv$ expresses on the following form, see the book~\cite[page~34]{Leonov94} where such a general model is introduced. This model also corresponds to a special case to the Oldroyd $8$-constant model, see~\cite{Oldroyd50}:
\begin{equation}\label{PEC5}
\overset{\bigtriangledown}{\bsigma} + a(\tr \bsigma) \bsigma + \eps (Dv : \bsigma) \bsigma = c(\tr \bsigma) \bdelta.
\end{equation}

\subsection{Mathematical results for the PEC model}

One of the key point in order to get the theoretical results is that a solution~$\bsigma$ of~\eqref{PEC5} can be express with respect to the conformation tensor, denoted by~$\bC$ (see the subsection~\ref{subsec:PEC-formulation} for more explanations):
\begin{equation*}\label{PEC6}
\bsigma = \gamma(\tr \bC) \bC
\quad \text{where} \quad  \overset{\bigtriangledown}{\bC} + \alpha(\tr \bC) \bC = \beta(\tr \bC) \bdelta,
\end{equation*}
the relations between the functions~$a$, $c$, the constant~$\eps$ and the functions $\alpha$, $\beta$, $\gamma$ will be given by the system~\eqref{assum4}.

\begin{example}\label{ex:2332}
For the specific case of the PEC model given by~\eqref{PEC3}, the functions~$a$, $c$ and the constant~$\eps$ are
$a\equiv 0$, $c\equiv 0$ and $\eps = \frac{2\xi'}{3G}$.
In that configuration, the corresponding functions $\alpha$, $\beta$ and $\gamma$ are
$\alpha\equiv 0$, $\beta\equiv 0$ and $\gamma(s) = \frac{3G}{3(1-\xi')+s\xi'}$.
We can note that such last relation appears in the works of Larson, see~\cite[equation (42), page 554]{Larson84}.
\end{example}
We are then first interested in the following type problem coupling the Navier-Stokes equations and the previously introduced constitutive law:
\begin{equation}\label{PEC8}
\left\{
\begin{aligned}
& \partial_t v + v\cdot\nabla v - \Delta v + \nabla p = \div \bsigma \\
& \div v = 0 \\[-0.2cm]
& \bsigma = \gamma(\tr \bC) \bC
\quad \text{where} \quad  \overset{\bigtriangledown}{\bC} + \alpha(\tr \bC) \bC = \beta(\tr \bC) \bdelta \\
& v\big|_{t=0} = v_{\mathrm{init}}, \qquad \bC\big|_{t=0} = \bC_{\mathrm{init}}.
\end{aligned}
\right.
\end{equation}

The first result concerns the uniqueness of solution to the above problem~\eqref{PEC8}. In a rigorous way we have (see Subsection~\ref{subsec:PEC-uniqueness} for the proof, and Appendix on page~\pageref{sec:notations} for the notations of the functional spaces):
\begin{proposition}[uniqueness and stability]\label{th:uniqueness}
Let $T>0$ and $d\in \{2,3\}$.
\par\noindent
We assume that the scalar functions~$\alpha$, $\beta$ and~$\gamma$ are of class~$\mathcal C^1$ on~$\mathbb R_+$.
\par\noindent
If $v_{\mathrm{init}}\in L^2$ with $\div v_{\mathrm{init}}=0$ in the sense of distributions, $\bC_{\mathrm{init}}\in L^\infty$ and if the system~\eqref{PEC8} possesses two solutions $(v_1,\bC_1)$ and $(v_2,\bC_2)$, such that, for $i\in \{1,2\}$,
\begin{equation*}
\begin{aligned}
& v_i \in L^\infty(0,T;L^2), && \quad \nabla v_i \in L^1(0,T;L^\infty), \\
& \bC_i \in L^\infty(0,T;L^\infty), && \quad \nabla \bC_i \in L^2(0,T;L^q) \quad \text{(for some $q>d$)}, 
\end{aligned}
\end{equation*}
then they coincide.\\
Moreover, the possible solution depends continuously on the initial conditions in the following sense: if $(v_i,\bC_i)$ is solution corresponding to the initial conditions $(v_{i,\mathrm{init}}, \bC_{i,\mathrm{init}})$, $i=\{1,2\}$, then there exists a constant $C\geq 0$ such that, for all $t\in (0,T)$,
$$
\|v_2(t)-v_1(t)\|_{L^2} + \|\bC_2(t)-\bC_1(t)\|_{L^2} \leq C \Big(  \|v_{2,\mathrm{init}} - v_{1,\mathrm{init}}\|_{L^2} + \|\bC_{2,\mathrm{init}}-\bC_{1,\mathrm{init}}\|_{L^2} \Big).
$$
\end{proposition}
The second result is an existence result of a local (in time) strong solution. It will be proved in Subsection~\ref{subsec:PEC-local-existence}:
\begin{proposition}[local existence]\label{th:local-existence}
Let $d\in \{2,3\}$, $r\in ]1,+\infty[$ and $q\in]d,+\infty[$.
\par\noindent
We assume that the scalar functions~$\alpha$, $\beta$ and~$\gamma$ are of class~$\mathcal C^1$ on~$\mathbb R_+$, and that there exists $k\geq 0$ such that for any $s>0$ we have
\begin{equation}\label{assum10}
|\alpha(s)| \lesssim 1+s^k, \qquad |\alpha'(s)| \lesssim 1+s^{k-1}, \qquad |\beta(s)| \lesssim 1+s^k \quad \text{and} \quad |\beta'(s)| \lesssim 1+s^{k-1}.
\end{equation}
If $v_{\text{init}} \in D^r_q$ and~$\bC_{\text{init}}\in W^{1,q}$ a symmetric positive definite almost everywhere then there exists $T_\star>0$ and a solution $(v,\bC)$ to the system~\eqref{PEC8} on~$[0,T_\star]$ satisfying:
\begin{equation*}
\begin{array}{ll}
 v \in L^r(0,T_\star;W^{2,q}),
& \partial_t v \in L^r(0,T_\star;L^q),\\
 \bC \in L^\infty(0,T_\star;W^{1,q}),
& \partial_t \bC \in L^r(0,T_\star;L^q).
\end{array}
\end{equation*}
\end{proposition}

\begin{remark}~\\[-0.6cm]
\begin{enumerate}
\item The regularity announced in Proposition~\ref{th:local-existence} allows us to apply Proposition~\ref{th:uniqueness}.
Indeed, if $v \in L^r(0,T_\star;W^{2,q})$ and $\partial_t v \in L^r(0,T_\star;L^q)$ then $v \in W^{1,r}(0,T_\star;L^q) \subset L^\infty(0,T_\star;L^2)$. The last inclusion is a consequence of the one-dimensional Sobolev embedding $W^{1,r}(0,T_\star)\subset L^\infty(0,T_\star)$ for $r>1$, and of the assumption $q>2$.
On the other hand, if $v \in L^r(0,T_\star;W^{2,q})$ then $\nabla v \in L^r(0,T_\star;W^{1,q})$. The Sobolev embedding $W^{1,q}\subset L^\infty$ right as soon as $q>d$, and the fact that $r>1$ implies that $\nabla v \in L^1(0,T_\star;L^\infty)$.
In a same way, if $\bC \in L^\infty(0,T_\star;W^{1,q})$ and $\partial_t \bC \in L^r(0,T_\star;L^q)$ then it is not difficult to prove that we have $\bC \in L^\infty(0,T_\star;L^\infty)$ and $\nabla \bC \in L^2(0,T_\star;L^q)$.
Consequently, the solution obtained in Proposition~\ref{th:local-existence} is unique.
\item As usually, it is certainly possible to obtain a global existence result assuming that the data are small enough. A proof, similar to that proposed in the subsection~\ref{subsec:PEC-local-existence}, could be based on the fact that for small data (that is, for small~$R_1$ and~$R_2$ - see the proof), the compact set~$K_{T,R_1,R_2}$ is stable by the function~$\Phi$ for all time $T>0$.
\item Similar results must be hold introducing a exterior force~$f$ in the right hand side of the first equation of~\eqref{PEC8} as soon as $f\in L_{\mathrm{loc}}^r(0,+\infty;L^q)$ (see for instance~\cite{chupin13} for similar result).
\end{enumerate}
\end{remark}

Next, we announce a global existence result if we add a few more assumptions (theses assumptions are satisfied for the models previously presented): 
\begin{proposition}[global existence]\label{th:global-existence}
Under the same assumptions that in Proposition~\ref{th:local-existence}, and if we assume that the reals $d$, $q$ and~$r$ satisfy
\begin{itemize}
\item $\displaystyle \frac{d}{2q}+\frac{1}{r}<\frac{1}{2}$,
\item $\displaystyle d=2 \quad \text{OR} \quad \big( d=3 \text{ but without taking into account the convective term $v\cdot \nabla v$} \big)$,
\end{itemize}
and that the functions $\alpha$, $\beta$ and $\gamma$ satisfy
\begin{itemize}
\item $\displaystyle \beta > 0 \quad \text{OR} \quad \beta \equiv 0$,
\item $\displaystyle \gamma>0 \quad \text{and} \quad s \in \mathbb R_+ \longmapsto s\gamma(s) \quad \text{is one to one bounded}$,
\item $\displaystyle \lim_{s\to + \infty} \Big( \frac{d\beta(s)}{s} - \alpha(s) \Big) < +\infty$,
\end{itemize}
then the solution to the system~\eqref{PEC8} given by Proposition~\ref{th:local-existence} holds taking any  $T_\star >0$.
\end{proposition}

Finally, we prove that the system~\eqref{PEC8} is equivalent to the following system, including the PEC model~\eqref{PEC3}:
\begin{equation}\label{PEC7}
\left\{
\begin{aligned}
& \partial_t v + v\cdot\nabla v - \Delta v + \nabla p = \div \bsigma \\
& \div v = 0 \\
& \overset{\bigtriangledown}{\bsigma} + a(\tr \bsigma) \bsigma + \eps (Dv : \bsigma) \bsigma = c(\tr \bsigma) \bdelta \\
& v\big|_{t=0} = v_{\mathrm{init}}, \qquad \bsigma\big|_{t=0} = \bsigma_{\mathrm{init}}.
\end{aligned}
\right.
\end{equation}
and can be conclude this section with the following main result:
\begin{theorem}[global existence - stress formulation]\label{th:global-existence-1}
Let $d\in \{2,3\}$, $r\in ]1,+\infty[$ and $q\in]d,+\infty[$ such that
\begin{itemize}
\item $\displaystyle \frac{d}{2q}+\frac{1}{r}<\frac{1}{2}$,
\item $\displaystyle d=2 \quad \text{OR} \quad \big( d=3 \text{ but without taking into account the convective term $v\cdot \nabla v$} \big)$.
\end{itemize}
Let $\eps>0$ and two scalar functions~$a$ and~$c$ of class~$\mathcal C^1$ on~$\mathbb R_+$. We assume that
\begin{itemize}
\item $\displaystyle c>0 \quad \text{OR} \quad c \equiv 0$,
\item $\displaystyle \frac{2}{\eps} a \Big( \frac{2}{\eps} \Big) \geq d \, c\Big( \frac{2}{\eps} \Big)$.
\end{itemize}
For any symmetric positive definite tensor $\bsigma_{\mathrm{init}}\in W^{1,q}$ such that
\begin{itemize}
\item $\displaystyle \eps\, \tr \bsigma_{\mathrm{init}} < 2$,
\end{itemize}
the system~\eqref{PEC7} admits a unique, stable, global strong solution.
\end{theorem}

\begin{remark}
\begin{enumerate}
\item The assumptions are naturally verified in the examples presented above. In particular, the last hypothesis is always true as soon as $\eps$ is small enough, whereas the previous hypothesis will be valid as soon as $xa(x)\geq d\, c(x)$ for $x$ large.
\item It is also important to note that assumption $\displaystyle \eps\, \tr \bsigma_{\mathrm{init}} < 2$ is practically not a smallness assumption on the data.
It can read
\begin{itemize}
\item either "for any initial data, there is $\varepsilon_0>0$ such that for $\varepsilon\in]0,\varepsilon_0]$ there is a unique strong global solution";
\item or "taking $\eps = \frac{2\xi}{\tr \bsigma_{\mathrm{init}}}$ where $\xi \in (0,1)$, the system~\eqref{PEC7} always admits a solution". This case corresponds to the PEC model, see \cite[Pages 553-554]{Larson84}.
\end{itemize}
\item It should be noted, however, that the result is not proven when $\eps=0$, and that a simple limit process $\eps \to 0$ cannot remedy it. Indeed, the estimates obtained in the proof show that, at~$\eps$ fixed, the solution is defined for all time~$t$ but can behave as $\mathrm{exp} (\mathrm{exp} (t/\eps))$.
\end{enumerate}
\end{remark}

\subsection{Proofs of the PEC model results}\label{sec:PEC-proof}

\subsubsection{Link between stress formulation and conformation formulation}\label{subsec:PEC-formulation}

In this subsection, we show that the conformation formulation of the PEC model gives a solution to the stress formulation.
This is one of the key points of the results since it allows the mathematical results from both the "stress models" framework and the "conformation model" framework to be used.
More precisely, we prove the following result:
\begin{proposition}\label{lemma:1053}
Let $\alpha$, $\beta$ and~$\gamma$ be three scalar functions of class~$\mathcal C^1$ on~$\mathbb R_+$.\\
We assume that $\gamma$ is positive and that $s\in\mathbb R_+ \longmapsto s\gamma(s)$ is one to one.\\
A tensor $\bsigma$ defined by~$\bsigma = \gamma(\tr \bC) \bC$ where $\bC$ is solution of
\begin{equation*}\label{eq:C1}
\left\{
\begin{aligned}
& \overset{\bigtriangledown}{\bC} + \alpha(\tr \bC) \bC = \beta(\tr \bC) \bdelta \\
& \bC\big|_{t=0} = \bC_{\mathrm{init}},
\end{aligned}
\right.
\end{equation*}
satisfies
\begin{equation}\label{eq:sigma1}
\left\{
\begin{aligned}
& \overset{\bigtriangledown}{\bsigma} + a(\tr \bsigma) \bsigma + b(\tr \bsigma) (Dv : \bsigma) \bsigma = c(\tr \bsigma) \bdelta \\
& \bsigma\big|_{t=0} = \gamma(\tr \bC_{\mathrm{init}}) \bC_{\mathrm{init}}.
\end{aligned}
\right.
\end{equation}
The relations between the functions~$a$, $b$, $c$ and $\alpha$, $\beta$, $\gamma$ are given by, for any $s>0$,
\begin{equation}\label{assum4}
\begin{aligned}
& a(s\gamma(s)) = \Big(1+\frac{s\gamma'(s)}{\gamma(s)} \Big) \alpha (s) - \frac{d\gamma'(s)}{\gamma(s)} \beta(s), \\
& b(s\gamma(s)) = -\frac{2\gamma'(s)}{\gamma(s)^2},\\
& c(s\gamma(s)) = \gamma(s) \beta(s).
\end{aligned}
\end{equation}
\end{proposition}

\paragraph{Proof of Proposition~\ref{lemma:1053}}

We consider a tensor $\bC$, whose trace is denoted by $s=\tr \bC$, satisfying
\begin{equation}\label{PEC11}
\overset{\bigtriangledown}{\bC} + \alpha(s)\, \bC = \beta(s)\, \bdelta.
\end{equation}
Taking the trace of~\eqref{PEC11}, we have
\begin{equation}\label{PEC12}
d_t s - 2 Dv : \bC + s \, \alpha(s) = d \, \beta(s).
\end{equation}
Now, we introduce the tensor $\bsigma$ defined by $\bsigma = \gamma(s)\, \bC$. Due to the definition of the upper convected derivative~\eqref{derive}, we have
\begin{equation*}\label{PEC13}
\overset{\bigtriangledown}{\bsigma} = \gamma(s) \overset{\bigtriangledown}{\bC} + \gamma'(s)\, d_t s \, \bC.
\end{equation*}
The equations~\eqref{PEC11} and~\eqref{PEC12} imply
\begin{equation*}\label{PEC14}
\overset{\bigtriangledown}{\bsigma} = \gamma(s) \Big(  - \alpha(s)\, \bC + \beta(s)\, \bdelta \Big) + \gamma'(s) \Big( 2 Dv : \bC - s \, \alpha(s) + d \, \beta(s) \Big)\bC.
\end{equation*}
Since $\bC=\bsigma/\gamma(s)$, we deduce that the quantity~$\bsigma$ satisfies the following equation
\begin{equation*}\label{PEC15}
\overset{\bigtriangledown}{\bsigma} = - \alpha(s)\, \bsigma + \gamma(s) \beta(s)\, \bdelta + \gamma'(s) \Big( 2 Dv : \frac{\bsigma}{\gamma(s)} - s \, \alpha(s) + d \, \beta(s) \Big)\frac{\bsigma}{\gamma(s)}.
\end{equation*}
It corresponds to~\eqref{eq:sigma1} as soon as the relations~\eqref{assum4} hold, where we have remarked that $\tr \bsigma = s \gamma(s)$. We also note that the functions $a$, $b$ and~$c$ are well defined by~\eqref{assum4} if $s\in \mathbb R_+ \longmapsto s\gamma(s) \in \mathbb R_+$ is an injection.
\finpreuve

\begin{remark}\label{rem:1102}
If we replace the upper convected derivative~$\overset{\bigtriangledown}{\cdot}$ by the derivative~$\mathcal D_\xi$, $-1\leq \xi \leq \xi$ (see~\eqref{derive1} for the notation), then Proposition~\ref{lemma:1053} becomes
\begin{align*}
& \bsigma = \gamma(\tr \bC) \bC, \quad \mathcal D_\xi \bC + \alpha(\tr \bC) \bC = \beta(\tr \bC) \bdelta \\
&\hspace{3cm} \Downarrow \\
& \mathcal D_\xi \bsigma + a(\tr \bsigma) \bsigma + \xi \, b(\tr \bsigma) (Dv : \bsigma) \bsigma = c(\tr \bsigma) \bdelta
\end{align*}
where the relations~\eqref{assum4} between~$a$, $b$, $c$ and~$\alpha$, $\beta$, $\gamma$ remain unchanged.
\end{remark}

\subsubsection{Proof of the uniqueness and stability result}\label{subsec:PEC-uniqueness}

In this subsection, we prove Theorem~\ref{th:uniqueness} which provides the uniqueness and stability of the solution to the system
\begin{equation*}\label{PEC19}
\left\{
\begin{aligned}
& \partial_t v + v\cdot\nabla v - \Delta v + \nabla p = \div \bsigma \\
& \div v = 0 \\
& \bsigma = \gamma(\tr \bC) \bC \\
& \overset{\bigtriangledown}{\bC} + \alpha(\tr \bC) \bC = \beta(\tr \bC) \bdelta \\
& v\big|_{t=0} = v_{\mathrm{init}} \quad \text{and} \quad \bC\big|_{t=0} = \bC_{\mathrm{init}}.
\end{aligned}
\right.
\end{equation*}
As usual, to obtain uniqueness result, we take the difference of the two solutions indexed by~$1$ and~$2$.
The vector $v=v_1-v_2$, the scalar $p=p_1-p_2$ and the tensors $\bsigma = \bsigma_1-\bsigma_2$, $\bC = \bC_1-\bC_2$ satisfy the following system:
\begin{equation}\label{PEC20}
\left\{
\begin{aligned}
& \partial_t v + v\cdot\nabla v_1 + v_2\cdot\nabla v - \Delta v + \nabla p = \div \bsigma \\
& \div v = 0 \\
& \bsigma = \gamma_1 \bC + (\gamma_1-\gamma_2) \bC_2 \\
& \partial_t \bC + v\cdot\nabla \bC_1 + v_2\cdot\nabla \bC - \bC\cdot\nabla v_1 - \bC_2\cdot\nabla v \\
& \qquad - \transp{(\nabla v)}\cdot \bC_1 - \transp{(\nabla v_2)}\cdot \bC + \alpha_1 \bC + (\alpha_1-\alpha_2) \bC_2 = (\beta_1-\beta_2) \bdelta,
\end{aligned}
\right.
\end{equation}
together with initial conditions $v\big|_{t=0} = v_{1,\text{init}} - v_{2,\text{init}}$ and $\bC\big|_{t=0} = \bC_{1,\text{init}} - \bC_{2,\text{init}}$. For the sake of simplification, we have noted~$\gamma_i$, $\alpha_i$ and $\beta_i$ the values of~$\gamma(\tr \bC_i)$, $\alpha(\tr \bC_i)$ and~$\beta(\tr \bC_i)$ respectively.

Taking the scalar product of the first equation of~\eqref{PEC20} by $v$ in $L^2$, we obtain
\begin{equation}\label{PEC21}
\begin{aligned}
\frac{1}{2}d_t \big( \|v\|_{L^2}^2 \big) + \|\nabla v\|_{L^2}^2 
& = -\int (v\cdot\nabla v_1)\cdot v - \int \bsigma:\nabla v \\
& \leq \|\nabla v_1\|_{L^\infty}\|v\|_{L^2}^2 + \frac{1}{2}\|\nabla v\|_{L^2}^2 + \frac{1}{2}\|\bsigma\|_{L^2}^2.
\end{aligned}
\end{equation}

The difference~$\bsigma$ of the stresses is controlled as follows
$|\bsigma| \leq |\gamma_1| |\bC| + |\gamma_1-\gamma_2| |\bC_2|$.
Since the function~$\gamma$ is of class~$\mathcal C^1$, and since the functions~$\bC_i$ are assumed to be bounded, the values~$\gamma(\tr \bC_i)$ are bounded.
Below, we denote by~$\sup(f)$ the bounded quantity $\sup |f(\tr \bC_i)|$ for any function $f\in \mathcal C(\mathbb R, \mathbb R)$. The differences of kind $\gamma_1-\gamma_2$ are so controlled by $\sup (\gamma') |\tr \bC_1-\tr \bC_2| \leq  \sup (\gamma') |\bC|$. We deduce
\begin{equation}\label{PEC21c}
\|\bsigma\|_{L^2} \leq \Big( \sup (\gamma) + \|\bC_2\|_{L^2} \sup(\gamma') \Big) \|\bC\|_{L^2}.
\end{equation}

Now, taking the scalar product of the fourth equation of~\eqref{PEC20} by~$\bC$ in~$L^2$, we have
\begin{equation}\label{PEC22}
\begin{aligned}
\frac{1}{2}d_t \big( \|\bC\|_{L^2}^2 \big) 
& = -\int (v\cdot\nabla \bC_1):\bC + \int (\bC\cdot\nabla v_1):\bC + \int (\bC_2\cdot\nabla v):\bC + \int (\transp{(\nabla v)}\cdot \bC_1):\bC \\
& \qquad + \int (\transp{(\nabla v_2)}\cdot \bC):\bC - \int \alpha_1 |\bC|^2 - \int (\alpha_1-\alpha_2) \bC_2:\bC + \int (\beta_1-\beta_2) \tr \bsigma.
\end{aligned}
\end{equation}
The first term of the right member can be estimated as follows, using successively Hölder's inequality, Sobolev's injection $H^1 \hookrightarrow L^{\frac{2q}{q-2}}$, holds for any $q>d$, Poincaré's inequality and Young's inequality (the constant $c_P$ being a universal constant involved in injections):
\begin{equation*}
\begin{aligned}
\Big| \int (v\cdot\nabla \bC_1):\bC \Big| & \leq \|v\|_{L^{\frac{2q}{q-2}}} \|\nabla \bC_1\|_{L^q} \|\bC\|_{L^2}\\
& \leq c_P\|\nabla v\|_{L^2} \|\nabla \bC_1\|_{L^q} \|\bC\|_{L^2}\\
& \leq \frac{1}{4}\|\nabla v\|_{L^2}^2 + c_P^2\|\nabla \bC_1\|_{L^q}^2 \|\bC\|_{L^2}^2.
\end{aligned}
\end{equation*}
Since the functions~$\alpha$ and $\beta$ are of class~$\mathcal C^1$, the equality~\eqref{PEC22} implies
\begin{equation}\label{PEC23}
\begin{aligned}
\frac{1}{2}d_t \big( \|\bC\|_{L^2}^2 \big)
& \leq \frac{1}{4}\|\nabla v\|_{L^2}^2 + c_P^2\|\nabla \bC_1\|_{L^q}^2 \|\bC\|_{L^2}^2 + 
 \|\nabla v_1\|_{L^\infty} \|\bC\|_{L^2}^2 + \frac{1}{4}\|\nabla v\|_{L^2}^2 \\
& \qquad + 2\|\bC_2\|_{L^\infty}^2 \|\bC\|_{L^2}^2 + 2\|\bC_1\|_{L^\infty}^2 \|\bC\|_{L^2}^2 + \|\nabla v_2\|_{L^\infty} \|\bC\|_{L^2}^2 \\
& \qquad + \sup(\alpha) \|\bC\|_{L^2}^2 + \sup(\alpha')\|\bC_2\|_{L^\infty} \|\bC\|_{L^2}^2 + \sup(\beta')\|\bC\|_{L^2}^2.
\end{aligned}
\end{equation}

Adding the estimates~\eqref{PEC21}, \eqref{PEC21c} and~\eqref{PEC23}, we deduce an estimate on the following form
\begin{equation*}\label{PEC24}
d_t \big( \|v\|_{L^2}^2 + \|\bC\|_{L^2}^2 \big) \leq h(t) ( \|v\|_{L^2}^2 + \|\bC\|_{L^2}^2),
\end{equation*}
where the function~$h$ is integrable on $(0,T)$. Indeed the function~$h$ contains quantities like $\|\nabla \bC_1\|_{L^q}^2$ or $\|\nabla v_1\|_{L^\infty} \|\bC_1\|_{L^\infty}^2$ which are, by the regularity assumption, integrable. From the initial condition we deduce that, for all $t\in (0,T)$
$$
\|v(t)\|_{L^2}^2 + \|\bC(t)\|_{L^2}^2 \leq \Big( \|v_{1,\text{init}} - v_{2,\text{init}}\|_{L^2}^2 + \|\bC_{1,\text{init}} - \bC_{2,\text{init}}\|_{L^2}^2 \Big) \mathrm{exp}\Big(\int_0^T h \Big),
$$
corresponding to the stability result. In particular, the uniqueness result also follows taking the same initial condition for the two solutions, i.e. $v_{\text{init}}=0$ and $\bC_{\text{init}} = \bzero$.
\finpreuve

\begin{remark}
Note that the uniqueness and stability results also holds if we replace the conditions 
\[\nabla v_i \in L^1(0,T;L^\infty) \quad \text{and} \quad \bC_i \in L^\infty(0,T;L^\infty)\]\\[-0.8cm]
by the conditions
\[\nabla v_i \in L^2(0,T;L^\infty) \quad \text{and} \quad \bC_i \in L^4(0,T;L^\infty).\]
The function~$h$ introduced at the end of the proof remains integrable on $(0,T)$.
\end{remark}

\subsubsection{Proof of the local existence result}\label{subsec:PEC-local-existence}

This subsection is devoted to the proof of the local existence result (Proposition~\ref{th:local-existence}) for the system~\eqref{PEC8} describing the fluid in term of velocity-conformation tensor $(v,\bC)$.
We rewrite the system~\eqref{PEC8} as a fixed point equation and apply Schauder’s theorem. More precisely, for given velocity vector~$\overline v$ and given tensor~$\overline \bC$, we are interested in the following system
\begin{equation}\label{PEC25}
\left\{
\begin{aligned}
& \partial_t v + \overline v\cdot\nabla \overline v - \Delta v + \nabla p = \div \bsigma \\
& \div v = 0 \\
& \bsigma = \gamma(\tr \overline \bC) \overline\bC \\
& \partial_t \bC + \overline v \cdot \nabla \bC - \bC\cdot \nabla \overline v - \transp{(\nabla \overline v)}\cdot \bC + \alpha(\tr \overline \bC) \bC = \beta(\tr \overline \bC) \bdelta.
\end{aligned}
\right.
\end{equation}
From classical results on the Stokes problem, see for instance~\cite{Giga91}, the velocity~$v$, solution of the two first equations of~\eqref{PEC25} satisfies:
\begin{equation}\label{PEC26}
\|v\|_{L^r(0,T;W^{2,q})} + \|\partial_t v\|_{L^r(0,T;L^q)} \lesssim \|v_{\text{init}}\|_{W^{2,q}} + \|\div \bsigma - \overline v \cdot \nabla \overline v\|_{L^r(0,T;L^q)}.
\end{equation}
We control~$\bsigma$ from~$\overline \bC$ using the third equation of~\eqref{PEC25} and introducing the continuous and non-decreasing real functions $\displaystyle \Gamma(s)=\max_{0<\tau<s} \gamma(\tau)$ and~$\displaystyle \widetilde \Gamma(s)=\max_{0<\tau<s} \gamma'(\tau)$:
\begin{equation*}\label{PEC26b0}
\begin{aligned}
\|\div \bsigma \|_{L^r(0,T;L^q)} \leq & \Gamma \big( \|\tr \overline \bC \|_{L^\infty(0,T;L^\infty)} \big) \|\overline \bC \|_{L^r(0,T;W^{1,q})} \\
& \qquad + \widetilde\Gamma \big( \|\tr \overline \bC \|_{L^\infty(0,T;L^\infty)} \big) \|\overline \bC \|_{L^\infty(0,T;L^\infty)} \|\nabla (\tr \overline \bC) \|_{L^r(0,T;L^q)}.
\end{aligned}
\end{equation*}
As a direct consequence of the definitions of trace and norms on tensorial spaces (see Appendix, page~\pageref{sec:notations}), we note that $|\tr \overline \bC |^2 \leq d |\overline \bC|^2$ and that $|\nabla(\tr \overline \bC)|^2 \leq d |\nabla \overline \bC|^2$. We deduce the following estimate
\begin{equation}\label{PEC26b}
\begin{aligned}
\|\div \bsigma \|_{L^r(0,T;L^q)} \leq & \Gamma \big( \sqrt{d} \|\overline \bC \|_{L^\infty(0,T;L^\infty)} \big) \|\overline \bC \|_{L^r(0,T;W^{1,q})} \\
& \qquad + \sqrt{d} \, \widetilde\Gamma \big( \sqrt{d} \|\overline \bC \|_{L^\infty(0,T;L^\infty)} \big) \|\overline \bC \|_{L^\infty(0,T;L^\infty)} \|\overline \bC \|_{L^r(0,T;W^{1,q})}.
\end{aligned}
\end{equation}
Given $\overline v$ and $\overline \bC$, the existence of a unique solution~$\bC$ to the fourth equation of~\eqref{PEC25} follows from the application of the method of characteristics. Some details are given in \cite[Appendix p. 26]{Fernandez-Guillen-Ortega} but we reformulate the estimates taking into account the specificities of the model.
We first compute the $L^2$ scalar product of the fourth equation of~\eqref{PEC25} by $q|\bC|^{q-2} \bC$. Then we take gradients in the fourth equation of~\eqref{PEC25} and compute the scalar product of the resulting equation with $q|\nabla \bC|^{q-2} \nabla \bC$. By addition, we find:
\begin{equation}\label{PEC27}
\begin{aligned}
d_t \big( \|\bC\|_{W^{1,q}}^q \big) \lesssim & \|\nabla \overline v\|_{L^\infty} \|\bC\|_{W^{1,q}}^q + \|\nabla^2 \overline v\|_{L^q} \|\bC\|_{L^\infty} \|\nabla \bC\|_{L^q}^{q-1} \\
& \qquad + (1+\|\overline \bC\|_{L^\infty}^k) \|\bC\|_{L^q}^{q-1} + (1+\|\overline \bC\|_{L^\infty}^{k-1}) \|\nabla \overline \bC\|_{L^q} \|\bC\|_{L^\infty} \|\nabla \bC\|_{L^q}^{q-1} \\
& \qquad  + (1+\|\overline \bC\|_{L^\infty}^k) \|\nabla \bC\|_{L^q}^q + (1+\|\overline \bC\|_{L^\infty}^{k-1}) \|\nabla \overline \bC\|_{L^q} \|\nabla \bC\|_{L^q}^{q-1}.
\end{aligned}
\end{equation}
We note that we used the assumptions $|\alpha(s)|\lesssim 1+s^k$, $|\beta(s)| \lesssim 1+s^{k}$, $|\alpha'(s)| \lesssim 1+s^{k-1}$ and $|\beta'(s)| \lesssim 1+s^{k-1}$.
Using Young inequalities as well as the Sobolev continuous embedding $W^{1,q} \hookrightarrow L^\infty$, valid for $q>d$, the estimate~\eqref{PEC27} now writes
\begin{equation}\label{PEC28}
d_t \big( \|\bC\|_{W^{1,q}}^q \big) \lesssim \big( 1 + \|\nabla \overline v\|_{L^\infty} + \|\nabla^2 \overline v\|_{L^q} + \|\overline \bC\|_{W^{1,q}}^k \big) \|\bC\|_{W^{1,q}}^q + 1 + \|\overline \bC\|_{W^{1,q}}^k.
\end{equation}
On the other hand, we express the time derivative~$\partial_t \bC$ from the fourth equation of~\eqref{PEC25} to deduce
\begin{equation}\label{PEC29}
\|\partial_t \bC\|_{L^q} \lesssim \big( 1 + \|\overline v\|_{W^{1,q}} + \|\overline C\|_{L^\infty}^k \big) \|\bC\|_{W^{1,q}} + 1 + \|\overline C\|_{L^\infty}^k.
\end{equation}
From~\eqref{PEC28} and~\eqref{PEC29}, we deduce that there exists a continuous function $\mathcal F$ such that
\begin{equation}\label{PEC30}
\|\bC\|_{L^\infty(0,T;W^{1,q})} + \|\partial_t \bC\|_{L^r(0,T;L^q)} \leq \mathcal F \big(\|\overline v\|_{L^1(0,T;W^{2,q})},\|\overline v\|_{L^r(0,T;W^{1,q})}, \|\overline \bC\|_{L^{rk}(0,T;W^{1,q})} \big).
\end{equation}

We conclude the proof introducing a compact convex set $K_{T,R_1,R_2}$ for which the application $\Phi : (\overline v,\overline \bC) \longmapsto (v,\bC)$ satisfies $\Phi(K_{T,R_1,R_2})\subset K_{T,R_1,R_2}$.
Roughly speaking, the set $K_{T,R_1,R_2}$ is composed of couples of functions $(w,\bB)$ such that
$$\|w\|_{L^r(0,T;W^{2,q})} + \|\partial_t w\|_{L^r(0,T;L^q)} \leq R_1,\quad \|\bB\|_{L^\infty(0,T;W^{1,q})} + \|\partial_t \bB\|_{L^r(0,T;L^q)} \leq R_2.$$
From the previous estimates~\eqref{PEC26}, \eqref{PEC26b} and~\eqref{PEC30}, if $(\overline v,\overline \bC) \in K_{T,R_1,R_2}$ then, for $R_1$ and $R_2$ large enough, and for~$T$ small enough then $(v,\bC)\in K_{T,R_1,R_2}$. This is the key to prove that the application~$\Phi$ has a fixed point, and then that the system~\eqref{PEC8} has a solution local in time, see~\cite{chupin13,Fernandez-Guillen-Ortega,Guillope-Saut-CRAS,Guillope-Saut3} for similar proofs.

\subsubsection{Proof of the global existence result}\label{subsec:PEC-global-existence}

In this subsection, we prove Proposition~\ref{th:global-existence}.
The proof consists in developing estimates obtained in Proposition~\ref{th:local-existence} for the velocity~$v$ and for the conformation tensor~$\bC$ which are finite for any finite time. It is decomposed into several parts.
The three first lemmas treat the unknowns velocity/stress/conformation as decoupled unknowns: Lemma~\ref{lem:conformation} gives estimates on the conformation~$\bC$. In Lemma~\ref{lem:stress} we obtain a bound for the stress whereas in Lemma~\ref{lem:velocity} we obtain regularity estimates for the velocity.
We conclude the proof of Proposition~\ref{th:global-existence} obtaining a bound on the solution~$(v,\bC)$ of the coupled system~\eqref{PEC8}, see Proposition~\ref{prop:bsigma}.

\begin{lemma}\label{lem:conformation}
Let $d\in \{2,3\}$, $q>d$ and $r>1$.

Let $\bC_{\mathrm{init}} \in W^{1,q}$ and~$v$ such that $\nabla v \in L^1(0,T_\star;L^\infty) \cap L^r(0,T_\star;W^{1,q})$.

Let~$\alpha$ and $\beta$ be two scalar functions of class $\mathcal C^1$ on $\mathbb R_+$ such that $\beta>0$ or $\beta \equiv 0$.

If the initial condition~$\bC_{\mathrm{init}}$ is a symmetric positive definite matrix, then the solution $\bC$ of the following equation
\begin{equation}\label{eq:conformation}
\left\{
\begin{aligned}
& \overset{\bigtriangledown}{\bC} + \alpha(\tr \bC) \bC = \beta(\tr \bC) \bdelta \\
& \bC\big|_{t=0} = \bC_{\mathrm{init}},
\end{aligned}
\right.
\end{equation}
remains a symmetric positive definite matrix for all time $t\in [0,T_\star]$ and all $x\in \T^d$.

Moreover, assuming
$$\displaystyle \lim_{s\to + \infty} \Big( \frac{d\beta(s)}{s} - \alpha(s) \Big) < +\infty,$$
we have, for any $t\in (0,T_\star)$:
\begin{equation}\label{PEC550}
\begin{aligned}
& \|\bC\|_{L^\infty(0,t;L^\infty)} \lesssim \mathrm{exp} \Big( \|\nabla v\|_{L^1(0,t;L^\infty)} \Big),\\
& \|\nabla \bC\|_{L^r(0,t;L^q)} \lesssim \Big( 1 + \|\nabla^2 v\|_{L^r(0,t;L^q)} \Big) \, \mathrm{exp} \Big( \|\nabla v\|_{L^1(0,T_\star;L^\infty)} \Big).
\end{aligned}
\end{equation}
\end{lemma}

{\bf Proof of Lemma~\ref{lem:conformation} -- }
The proof of the positiveness of the conformation tensor~$\bC$ is given by M.A. Hulsen (see~\cite{Hulsen90}) in the three dimensional case. It can adapt as it is in the two dimensional case.

In order to obtain the first estimate of~\eqref{PEC550}, we take the scalar product of~\eqref{eq:conformation} by $q|\bC|^{q-2} \bC$:
\begin{equation*}
d_t|\bC|^q + q \alpha(\tr \bC) |\bC|^q = q \beta(\tr \bC) \tr \bC |\bC|^{q-2} + q|\bC|^{q-2}(\bC\cdot \nabla v + \transp{\nabla v}\cdot \bC) : \bC.
\end{equation*}
Since $\bC$ is symmetric positive definite, we can use the equivalence between~$\tr \bC$ and~$|\bC|$, see relations~\ref{eq:equivalent}. We deduce
\begin{equation}\label{PEC51}
d_t|\bC|^q \leq 2q |\nabla v| |\bC|^q + q \Big( \frac{d\beta(\tr \bC)}{\tr \bC} - \alpha(\tr \bC) \Big) |\bC|^q.
\end{equation}
If $\displaystyle \lim_{s\to + \infty} \Big( \frac{d\beta(s)}{s} - \alpha(s) \Big) < +\infty$ then, due to the positiveness of $\tr \bC$ and the continuity of~$\alpha$ and~$\beta$, we have
\begin{equation*}
\frac{d\beta(\tr \bC)}{\tr \bC} - \alpha(\tr \bC) \lesssim 1.
\end{equation*}
The estimate~\eqref{PEC51} becomes
\begin{equation*}\label{PEC52}
d_t|\bC|^q \lesssim 2q |\nabla v| |\bC|^q + q |\bC|^q.
\end{equation*}
Note that the symbol $\lesssim$ introduced here does not depend on the integer~$q$. Integrating with respect the the space variable, we have
\begin{equation*}\label{PEC53}
d_t\|\bC\|_{L^q}^q \lesssim q \big( 2 \|\nabla v\|_{L^\infty} + 1 \big) \|\bC\|_{L^q}^q.
\end{equation*}
We then have $d_t\|\bC\|_{L^q} \lesssim \big( 2 \|\nabla v\|_{L^\infty} + 1 \big) \|\bC\|_{L^q}$ so that
\begin{equation*}\label{PEC54}
\|\bC\|_{L^\infty(0,T;L^q)} \lesssim \|\bC_{\text{init}}\|_{L^q} \mathrm{exp} \Big( \int_0^T \big( 2\|\nabla v\|_{L^\infty} + 1 \big) \Big) \lesssim \mathrm{exp} \Big( \|\nabla v\|_{L^1(0,t;L^\infty)} \Big).
\end{equation*}
Since $\bC_{\text{init}} \in W^{1,q} \subset L^\infty$, the last estimate does not depend on $q$, passing to the limit $q\to+ \infty$, we deduce the first estimate.\\

Finally, to obtain the second estimate of~\eqref{PEC550}, we now take the gradients of~\eqref{eq:conformation} and compute the scalar product in $\mathbb R^{d\times d \times d}$ with $q |\nabla \bC|^{q-2} \nabla \bC$:
\begin{equation*}\label{PEC56}
\begin{aligned}
d_t |\nabla \bC|^q \leq & 3q |\nabla v| |\nabla \bC|^q + 2q |\nabla ^2 v| |\bC| |\nabla \bC|^{q-1} + q |\alpha(\tr \bC)| |\nabla \bC|^q  \\
&\qquad + q |\alpha'(\tr \bC)| |\nabla (\tr \bC)| |\bC| |\nabla \bC|^{q-1} + q \sqrt{d} |\beta'(\tr \bC)| |\nabla (\tr \bC)| | \nabla \bC|^{q-1}.
\end{aligned}
\end{equation*}
Note that from the previous bound $\|\bC\|_{L^\infty(0,T;L^\infty)} \lesssim 1$, we have $|f(\tr \bC)| \lesssim 1$ for any continuous function~$f$, like $\alpha$, $\alpha'$ or~$\beta'$.
We deduce the following estimate
\begin{equation*}\label{PEC57}
d_t |\nabla \bC|^q \lesssim (|\nabla v| +1) |\nabla \bC|^q + |\nabla ^2 v| |\nabla \bC|^{q-1}.
\end{equation*}
Here, the coefficient hidden behind the symbol $\lesssim$ depends on the integer~$q$ (but that does not matter here).
After integration in the space variable, we deduce
\begin{equation*}\label{PEC58}
d_t \|\nabla \bC\|_{L^q}^q \lesssim (\|\nabla v\|_{L^\infty} +1) \|\nabla \bC\|_{L^q}^q + \|\nabla ^2 v\|_{L^q} \|\nabla \bC\|_{L^q}^{q-1}.
\end{equation*}
We multiply by $\|\nabla \bC\|_{L^q}^{r-q}$ and use the Young inequality $ab^{r-1}\leq \frac{1}{r}a^r + \frac{r-1}{r}b^r$ in order to have
\begin{equation*}\label{PEC59}
d_t \|\nabla \bC\|_{L^q}^r \lesssim (\|\nabla v\|_{L^\infty} +1) \|\nabla \bC\|_{L^q}^r + \|\nabla ^2 v\|_{L^q}^r.
\end{equation*}
The Gronwall lemma allows us to conclude.
\finpreuve

\begin{remark}
Note also that the proof of the positiveness of the conformation tensor~$\bC$ presented by Hulsen use the assumption $\beta>0$. But it is not difficult to adapt the proof to the case where~$\beta$ is identically zero. On the other hand, if~$\beta$ changes sign, the result can be false. We will see in Section~\ref{sec:conclusion}, page~\pageref{sec:false}, an example where the term~$\beta$ changes sign, and where~$\bC$ does not remain positive.
\end{remark}

\begin{lemma}\label{lem:stress}
Let $\bsigma:(0,T_\star)\times \T^d \longrightarrow \mathbb R^{d\times d}$ be defined by
\begin{equation*}\label{eq:stress}
\left\{
\begin{aligned}
& \bsigma = \gamma(\tr \bC) \bC \\
& \overset{\bigtriangledown}{\bC} + \alpha(\tr \bC) \bC = \beta(\tr \bC) \bdelta \\
& \bC\big|_{t=0} = \bC_{\mathrm{init}}.
\end{aligned}
\right.
\end{equation*}
We assume that~$\alpha$,  $\beta$ and $\gamma$ are three scalar functions of class $\mathcal C^1$ on $\mathbb R_+$ such that $\beta>0$ or $\beta \equiv 0$, and
\begin{equation*}\label{assum1bis}
s\longmapsto s\gamma (s) \quad \text{is bounded}.
\end{equation*}
If the initial condition~$\bC_{\mathrm{init}}\in W^{1,q}$ is a symmetric positive definite matrix then the stress~$\bsigma$ is bounded:
\begin{equation*}\label{PEC31}
\|\bsigma\|_{L^\infty(0,T_\star;L^\infty)} \lesssim 1.
\end{equation*}
\end{lemma}

{\bf Proof of Lemma~\ref{lem:stress} -- }
We perform the proof of the $L^\infty$-estimate of the stress~$\bsigma$ point by point (for simplicity we note here~$\bC$ and~$\bsigma$ instead of~$\bC(t,x)$ and~$\bsigma(t,x)$ respectively).
Due to first part of Lemma~\ref{lem:conformation}, we know that~$\bC$, and then~$\bsigma$ are symmetric positive definite matrices.
Consequently, $|\bsigma| \leq \tr \bsigma$ (see relations~\ref{eq:equivalent}) and since we have defined the stress by $\bsigma = \gamma(\tr \bC)\bC$, we deduce that $|\bsigma| \leq \tr \bC \, \gamma (\tr \bC)$.
The assumption on~$\gamma$ allows us to conclude that $|\bsigma| \lesssim 1$.
\finpreuve

\begin{lemma}\label{lem:velocity}
Let $d\in \{2,3\}$, $r\in ]1,+\infty[$ and $q\in]d,+\infty[$ such that $\frac{d}{2q}+\frac{1}{r}<\frac{1}{2}$.
We assume that
\begin{equation*}\label{assum3bis}
d=2 \quad \text{OR} \quad \Big( d=3 \text{ but without taking into account the convective term $v\cdot \nabla v$} \Big).
\end{equation*}
If $v_{\mathrm{init}} \in D_q^r$ and $\bsigma \in L^r(0,T_\star;W^{1,q}) \cap L^\infty(0,T_\star;L^\infty)$
then the solution~$v$ of the following Navier-Stokes equation
\begin{equation*}\label{PEC32}
\left\{
\begin{aligned}
& \partial_t v + v\cdot\nabla v - \Delta v + \nabla p = \div \bsigma \\
& \div v = 0 \\
& v\big|_{t=0} = v_{\mathrm{init}},
\end{aligned}
\right.
\end{equation*}
satisfies the estimates, for any $t\in (0,T_\star)$:
\begin{equation}\label{PEC33}
\begin{aligned}
& \|\nabla v\|_{L^\infty(0,t;L^\infty)} \lesssim 1 + \|\bsigma\|_{L^\infty(0,t;L^\infty)} \ln (\mathrm e + \|\nabla \bsigma\|_{L^r(0,t;L^q)}),\\
& \|\nabla^2 v\|_{L^r(0,t;L^q)} \lesssim \|\nabla \bsigma\|_{L^r(0,t;L^q)}.
\end{aligned}
\end{equation}
\end{lemma}

{\bf Proof of Lemma~\ref{lem:velocity} -- }
The proof is based on the integral representation of the solution to the Navier-Stokes equation:
\begin{equation}\label{PEC34}
v(t,x) =\mathrm e^{t \Delta} v_{\text{init}} + \int_0^t \mathrm e^{(t-s) \Delta} \mathds{P} \div f(s,x)\, \mathrm ds,
\end{equation}
where we have introduced $f:=\bsigma - v\otimes v$.\\[0.2cm]
$\phantom{toto}\bullet$ The first estimate is obtain as follows: Due to the analyticity of the heat kernel, and using the Sobolev injections (essentially $W^{\frac{d}{q},q} \hookrightarrow L^\infty$), we successively prove that
\[
\|\mathrm e^{t\Delta} \Delta^\alpha g\|_{L^\infty} \lesssim t^{-\alpha} \|g\|_{L^\infty}
\quad \text{and} \quad
\|\mathrm e^{t\Delta} \Delta g\|_{L^\infty} \lesssim t^{-\frac{q+d}{2q}} \|\nabla g\|_{L^q}.
\]
These two estimates make it possible, from~\eqref{PEC34} and as soon as $\frac{d}{2q}+\frac{1}{r}<\frac{1}{2}$, to obtain (see~\cite{Chupin14,Constantin07} for similar results):
\begin{equation}\label{PEC35}
\|\nabla v(t,\cdot)\|_{L^\infty} \lesssim  1 + \| f \|_{L^\infty(0,t;L^\infty)} \ln \big( \mathrm e + \|\nabla f\|_{L^r(0,t;L^q)} \big).
\end{equation}
In the two dimensional case ($d=2$), an analysis of the Navier-Stokes equation implies that (see~\cite{Chupin14,Constantin07}) for a.e. $t\in (0,T_\star)$:
\begin{equation*}\label{PEC36}
\begin{aligned}
& \|v\|_{L^\infty(0,t;L^\infty)} \lesssim 1, \\
& \|\nabla v\|_{L^r(0,t;L^q(\T^2))} \lesssim 1 \quad \text{for any $1<q,r<+\infty$}. 
\end{aligned}
\end{equation*}
This makes it possible to control the non-linearity $v\otimes v$ and reduce the control of~$f$ (resp.~$\nabla f$) to that of~$\bsigma$ (resp.~$\nabla \bsigma$). It results the first estimate of~\eqref{PEC33}.\par
In the three dimensional case ($d=3$), we do not take into account the contribution $v\otimes v$ so that $f=\bsigma$ and the results directly follow from~\eqref{PEC35}.\\[0.2cm]
$\phantom{toto}\bullet$ The second estimate announced in Lemma~\ref{lem:velocity} comes from the integral representation~\eqref{PEC34} again and to the fact that the linear operator
$g \mapsto \int_0^t \mathrm e^{(t-s) \Delta} \Delta g(s) \, \mathrm ds$,
is bounded in $L^r(0,T;L^q)$ for $1<q,r<+\infty$, see~\cite[p. 64]{Lemarie}.
\finpreuve

\begin{proposition}\label{prop:bsigma}
Under the assumptions of Proposition~\ref{th:global-existence}, the solution  $(v,\bC)$ to the system~\eqref{PEC8} on~$[0,T_\star]$ satisfies, for any $t\in (0,T_\star)$,
\begin{equation*}\label{PEC330}
\begin{aligned}
& \|v\|_{L^r(0,t;W^{2,q})} \lesssim 1,\quad && \|\partial_t v\|_{L^r(0,t;L^q)} \lesssim 1, \\
& \|\bC\|_{L^\infty(0,t;W^{1,q})} \lesssim 1,\quad && \|\partial_t \bC\|_{L^r(0,t;L^q)} \lesssim 1.
\end{aligned}
\end{equation*}
\end{proposition}

{\bf Proof of Proposition~\ref{prop:bsigma} -- }
Combining the results of Lemma~\ref{lem:velocity} and Lemma~\ref{lem:stress}, we deduce that for any $t\in (0,T_\star)$ we have
\begin{equation}\label{PEC33bis}
\begin{aligned}
& \|\nabla v\|_{L^\infty(0,t;L^\infty)} \lesssim 1 + \ln (\mathrm e + \|\nabla \bsigma\|_{L^r(0,t;L^q)}),\\
& \|\nabla^2 v\|_{L^r(0,t;L^q)} \lesssim \|\nabla \bsigma\|_{L^r(0,t;L^q)}.
\end{aligned}
\end{equation}
The first goal is then to obtain a bound on
$$y(t) := \|\nabla \bsigma\|_{L^r(0,t;L^q)}^r.$$
We use the stress formulation introduced in Subsection~\ref{subsec:PEC-formulation}, see Proposition~\ref{lemma:1053}:
\begin{equation}\label{eq:1552}
\overset{\bigtriangledown}{\bsigma} + a(\tr \bsigma) \bsigma + b(\tr \bsigma) (Dv : \bsigma) \bsigma = c(\tr \bsigma) \bdelta.
\end{equation}
We first take gradients in equation~\eqref{eq:1552} and compute the~$L^2$ scalar product of the resulting equation with~$|\nabla \bsigma|^{q-2} \nabla \bsigma$.
From the bound $\|\bsigma\|_{L^\infty(0,T_\star;L^\infty)} \lesssim 1$, we control $f(\tr \bsigma)$ for $f=a,b,c,a',b',c'$.
Using the Hölder inequality, we find
\begin{equation*}\label{PEC37}
d_t \|\nabla \bsigma\|_{L^q}^q \lesssim (1 + \|\nabla v\|_{L^\infty}) \|\nabla \bsigma\|_{L^q}^q + \|\nabla^2 v\|_{L^q} \|\nabla \bsigma\|_{L^q}^{q-1}.
\end{equation*}
We next multiply this result by $\|\nabla \bsigma\|_{L^q}^{r-q}$ and use the Young type inequality $xy^{r-1}\leq \frac{1}{r}x^r + (1-\frac{1}{r})y^r$ to treat the last term. We deduce
\begin{equation*}\label{PEC38}
d_t \|\nabla \bsigma\|_{L^q}^r \lesssim (1 + \|\nabla v\|_{L^\infty}) \|\nabla \bsigma\|_{L^q}^r + \|\nabla^2 v\|_{L^q}^r.
\end{equation*}
Using estimate~\eqref{PEC33bis}, we conclude that
$y'(t) \lesssim y(t) \ln (\mathrm e + y(t)) + y(t)$
that implies that $y(t)$ is well defined for any time $t>0$ (it is bounded by a double exponential function of kind $\mathrm e^{\mathrm e^t}$). From Lemma~\ref{lem:velocity}, we then have for any $t\in [0,T_\star]$:
$$\|v\|_{L^r(0,t;W^{2,q})} \lesssim y(t) \lesssim 1.$$
We also note that, from~\eqref{PEC33bis}, we have
$$\|\nabla v\|_{L^\infty(0,t;L^\infty)} \lesssim 1.$$
In order to obtain the bound on the conformation~$\bC$, we use Lemma~\ref{lem:conformation} and the two previous estimates:
$$\|\bC\|_{L^\infty(0,t;W^{1,q})} \lesssim 1.$$
The two last bounds follow from the equations~\eqref{PEC8} expressing $\partial_t v$ and~$\partial_t \bC$ with respect to $v$, $\nabla v$, $\Delta v$, $\div \bsigma$, $\bC$ and~$\nabla \bC$.
\finpreuve

\subsubsection{Proof of the main theorem \ref{th:global-existence-1}}\label{subsec:PEC-proof}

From Proposition~\ref{th:global-existence} we know that there exists a strong solution to the conformation formulation~\eqref{PEC8} as soon as some assumptions on the functions~$\alpha$, $\beta$ and~$\gamma$ are satisfied.
Moreover, in Subsection~\ref{subsec:PEC-formulation}, we proved that any solution to the conformation formulation implies a solution to the stress formulation~\eqref{PEC7}, the link between functions~$\alpha$, $\beta$, $\gamma$ and functions $a$, $b$, $c$ being given by~\eqref{assum4}.
To prove the Theorem~\ref{th:global-existence-1} it suffices to show that the assumptions on functions~$a$ and~$c$ correspond to assumptions on $\alpha$, $\beta$ and~$\gamma$ in Proposition~\ref{th:global-existence}. Note that this approach does not show the uniqueness of the solution. However, the proof of uniqueness may be made independently by following exactly the same method as for proof of existence in the conformation formulation, see Subsection~\ref{subsec:PEC-uniqueness}.

\begin{itemize}
\item Let $a$, $c$ be two scalar functions of class~$\mathcal C^1$ on $\mathbb R_+$, and $\eps>0$.
We introduce the three scalar functions $\alpha$, $\beta$ and~$\gamma$ defined on $\mathbb R_+$ by
$$
\alpha(s) = (1+\eps\, s) \Big( a(s\gamma(s)) - \frac{d\, \eps}{2} c(s\gamma(s)) \Big),
\qquad
\beta(s) = \frac{c(s\gamma(s))}{\gamma(s)}
\quad \text{and} \quad
\gamma(s) = \frac{2}{1+\eps \, s},
$$
so that the relations~\eqref{assum4} hold with $b = \eps$.
\item Let $\bsigma_{\mathrm{init}}\in W^{1,q}$ be a symmetric positive definite tensor. Assuming $\eps\, \tr \bsigma_{\mathrm{init}} < 2$ it is possible to define the symmetric positive definite matrix~$\bC_{\mathrm{init}}\in W^{1,p}$ by
$$
\bC_{\mathrm{init}} = \frac{\bsigma_{\mathrm{init}}}{2-\eps\, \tr \bsigma_{\mathrm{init}}},
$$
so that $\bsigma_{\mathrm{init}} = \gamma(\tr \bC_{\mathrm{init}}) \bC_{\mathrm{init}}$.
\end{itemize}

Consequently, due to Subsection~\ref{subsec:PEC-formulation}, the solution to~\eqref{PEC8} will be solution to~\eqref{PEC7}.

We must therefore prove that functions $\alpha$, $\beta$ and~$\gamma$ defined above satisfy the assumptions of Proposition~\ref{th:global-existence}.

Since $a$ and~$c$ are functions of class~$\mathcal C^1$ on~$\mathbb R_+$, and due to the definition of~$\gamma$, it is clear that $\alpha$, $\beta$ and~$\gamma$ are of class~$\mathcal C^1$ too.
Moreover the assumption $\big( c>0 \quad \text{OR} \quad c \equiv 0 \big)$ clearly implies $\big( \beta>0 \quad \text{OR} \quad \beta \equiv 0 \big)$.

Finally, the other assumptions are direct consequences of the following asymptotic behaviors, for $s\to +\infty$:
$$
\gamma(s) = \frac{2}{\eps s} + \mathcal O\Big( \frac{1}{s^2} \Big), \quad 
\alpha(s) = \eps \Big( a \Big( \frac{2}{\eps} \Big) - \frac{d\, \eps}{2} c\Big( \frac{2}{\eps} \Big)  \Big) s + \mathcal O(1) \quad \text{and} \quad
\beta(s) = \frac{\eps}{2} c\Big( \frac{2}{\eps} \Big) s + \mathcal O(1).
$$
It should be noted in particular that hypothesis
$\displaystyle \lim_{s\to + \infty} \Big( \frac{d\beta(s)}{s} - \alpha(s) \Big) < +\infty$ is verified because we have assume that $\displaystyle \frac{2}{\eps} a \Big( \frac{2}{\eps} \Big) \geq d \, c\Big( \frac{2}{\eps} \Big)$.

\section{Global existence results for other models?}\label{sec:others}

In this section, we highlight the fact that the method proposed above to obtain the existence of a solution can be adapted to other non-Newtonian models. Some are direct corollaries of previous results (as is the case for PTT-type models, see Subsection~\ref{subsec:PTT}) while others require some intermediate results (for example to be able to apply the method to the MGI-type model, see Subsection~\ref{subsec:MGI} or pom-pom type polymer, Subsection~\ref{subsec:pom-pom}).
We end this section by showing almost the same type of models for which the positivity of the conformation is not respected over time, see Subsection~\ref{sec:false}.

\subsection{The Phan-Thien-Tanner (PTT) model with additive nonlinear term}\label{subsec:PTT}

The crucial point in modeling is the choice of a suitable constitutive equation which has a capability to correctly represent non-linear behavior of the melts. In recent years, significant progress has been made to develop such constitutive equations. Usual nonlinearities from the Oldroyd model, see equation~\eqref{eq:Oldroyd}, can be described by the PTT model.
Such model is an extension of the Oldroyd model to include a function dependent upon $\tr \bsigma$, the trace of the polymer stress.
More precisely, the constitutive law for the PTT model (with the additive term those the coefficient is $\widetilde \eps >0$) is
\begin{equation}\label{O10}
\lambda \overset{\bigtriangledown}{\bsigma} + f(\tr \bsigma) \bsigma + \widetilde \eps (Dv:\bsigma)(\lambda \bsigma + \eta_p \bdelta)= 2 \eta_p Dv.
\end{equation}
There are three forms of the function $f(\tr \bsigma)$ found in the literature:
\begin{equation*}
f(s) = \left\{
\begin{aligned}
& 1 + \frac{\kappa \lambda}{\eta_p}s &&\qquad \text{Linear PTT}\\
& 1 + \frac{\kappa \lambda}{\eta_p}s + \frac{1}{2}\Big( \frac{\kappa \lambda}{\eta_p}s \Big)^2 &&\qquad \text{Quadratic PTT}\\
& \mathrm{exp}\Big( \frac{\kappa \lambda}{\eta_p}s \big) &&\qquad \text{Exponetial PTT},
\end{aligned}
\right.
\end{equation*}
where $\kappa \in [0,1]$ is a model parameter.
Parameter $\kappa$ is inversely proportional to the extensional viscosity of the fluid and the linear or quadratic models only approach well the exponential form at low deformations.
The linear and exponential forms of the PTT model are extensively used, and are the two forms mentioned in~\cite[Sections 5.6.5 and 5.6.6]{Tanner00} about the PTT model. The exponential model was first proposed by Phan-Thien~\cite{PhanThien78} a year after the linear model of~\cite{PTT77}. The quadratic form is far less widely used or mentioned in the literature, but is used, for example, to model the wire-coating process in~\cite{NW02}, where all three PTT forms are investigated.

We remark that this type of model~\eqref{O10} falls within the framework that we have described in this article, and that consequently we can apply the result of Theorem~\ref{th:global-existence-1}.

Indeed, introducing $\widetilde \bsigma = \lambda \bsigma + \eta_p \bdelta$, the equation~\eqref{O10} takes the form of the third equation of~\eqref{PEC7} with
\[a(s) = \frac{1}{\lambda} f\Big( \frac{s-d\eta_p}{\lambda} \Big), \qquad c(s)= \frac{\eta_p}{\lambda} f\Big( \frac{s-d\eta_p}{\lambda} \Big) \quad \text{and} \quad \eps = \frac{\widetilde \eps}{\lambda}.\]
Under the condition $\widetilde \eps \, d \, \eta_p \leq 2\lambda$, corresponding to the condition $\frac{2}{\eps}a\big( \frac{2}{\eps} \big) \geq d\, c\big( \frac{2}{\eps} \big)$, it is then possible to apply Theorem~\ref{th:global-existence-1} and deduce the existence of an unique global strong solution.

\begin{remark}~\\[-0.5cm]
\begin{enumerate}
\item The same results are also true with other classic models:
\begin{itemize}
\item[i-] the FENE-P (Finitely Extensible Non-Linear Elastic-Peterlin) proposed by Bird, Armstrong and Hassager~\cite{Bird87}. See also recent theoretical results relative to a diffusive case~\cite{Renardy17};
\item[ii-] the FENE-CR ( Finitely Extensible Non-Linear Elastic-Chilcott and Rallison) proposed by Chilcott and Rallison~\cite{Chilcott88}.
\end{itemize}
If we add a non-linear term of the form $\varepsilon (Dv:\bsigma)\bsigma$ (with a possibly very small~$\varepsilon$ positive coefficient) then we have a global existence of a strong solution for any initial data.
\item In the previous models (PEC, MGI, PTT, FENE-P or FENE-CR) it is possible to use the derivative~$\mathcal D_\xi$, $\xi \neq 0$, instead of the UCM derivative (see~\eqref{derive1} for the notation). All the proof presented in this paper remains correct.
In particular, introducing $\widetilde \bsigma = \bsigma + \frac{\eta_p}{\xi \lambda} \bdelta$, the equation~\eqref{O10} - with~$\mathcal D_\xi \bsigma$ instead of~$\overset{\bigtriangledown}{\bsigma}$ - takes the form of equation~\eqref{PEC7}$_3$ where we replace $\varepsilon$ by $\xi\, \varepsilon$.
\end{enumerate}
\end{remark}

\subsection{The Marrucci, Greco and Ianniruberto (MGI) type model}\label{subsec:MGI}

The incorporation of additional mechanisms, such as contour length fluctuations and stress release, see~\cite{Doi88}, leads to a precise description of the viscoelastic properties under the flow. 
Many modifications such as taking into account the finite stretchability of the chain have also been introduced. Thus, in~\cite{Marrucci01}, Marrucci, Greco and Ianniruberto argued that a force balance on the entanglement nodes should be fulfilled, resulting in a modified~$\bQ$ tensor. In this manner, a better agreement with experimental data for the normal stress ratio could be obtained. The modified $\bQ$ tensor is given with respect to the Finger tensor~$\bC$, measuring the deformation of a fluid element, by
\begin{equation}\label{eq:1356}
\bQ = \frac{\sqrt{\bC}}{\tr{\sqrt{\bC}}} \qquad \text{with} \qquad \overset{\bigtriangledown}{\bC} = \bzero.
\end{equation}
The tensor $\sqrt{\bC}$ here appearing is a special case of the Seth tensor $\bC^a$ sometimes proposed in phenomenological equations for rubber or for viscoelastic liquids, see~\cite{Larson13}. Larson~\cite{Larson97} relates the $\sqrt{\bC}$ tensor to the disclination lines liquid crystal, which also carry a constant tension.
The other characteristic of the choice~\eqref{eq:1356} is the ratio defining~$\bQ$. This point is essential from a theoretical point of view since it clearly implies a $L^\infty$ bound on the tensor~$\bQ$. According to Marrucci, Greco and Ianniruberto~\cite{Marrucci00}, this form of writing is a generalization of the linear case "$\bQ=\bC$" to large deformations. More exactly, they claim that during retraction, the number of strands per unit volume decreases in inverse proportion to the average increase of strand length scale with~$\tr \sqrt{\bC}$.
The stress tensor~$\bsigma$ obeys time-strain separability and is obtain as
\[\bsigma = G \, \Lambda(t) \, \bQ,\]
where $G>0$ is a physical constant and $\Lambda$ is a function.
In practice, the function~$\Lambda$ is a solution of a differential equation coupled with the velocity and/or the stress, see for instance~\cite{Marrucci01,Verbeeten01,Verbeeten02}. Nevertheless, following the remark given in~\cite[equation (7), page 100]{Marrucci01}: {\it "The reptation function~$\Lambda$ is not very different from a single exponential $\Lambda(t) = \mathrm e^{-t/\tau_0}$"}, we will assume that the function~$\Lambda$ is any positive and regular function. We are therefore interested in the following system
\begin{equation}\label{MGI0}
\left\{
\begin{aligned}
& \partial_t v + v\cdot\nabla v - \Delta v + \nabla p = \div \bsigma \\
& \div v = 0 \\
& \bsigma = \Lambda(t) \frac{\sqrt{\bC}}{\tr \sqrt{\bC}} \quad \text{where} \quad \overset{\bigtriangledown}{\bC} = 0 \\
& v\big|_{t=0} = v_{\mathrm{init}}, \qquad \bC\big|_{t=0} = \bC_{\mathrm{init}}.
\end{aligned}
\right.
\end{equation}

\subsubsection{Mathematical results for the MGI model}

Following the same ideas as in the previous Section, we have the following results:
\begin{proposition}[uniqueness and stability]\label{th:uniqueness1}
Let $T>0$ and $d\in \{2,3\}$.\par\noindent
We assume that~$\Lambda$ is a positive function of class~$\mathcal C^1$ on~$\mathbb R_+$.\par\noindent
If $v_{\mathrm{init}}\in L^2$ with $\div v_{\mathrm{init}}=0$ in the sense of distributions, $\bC_{\mathrm{init}}\in L^2$ is symmetric positive definite almost everywhere, and if the system~\eqref{MGI0} possesses two solutions $(v_1,\bC_1)$ and $(v_2,\bC_2)$, such that, for $i\in \{1,2\}$,
\begin{equation*}
\begin{aligned}
& v_i \in L^\infty(0,T;L^2), && \quad \nabla v_i \in L^1(0,T;L^\infty), \\
& \bC_i \in L^\infty(0,T;L^\infty), && \quad \nabla \bC_i \in L^2(0,T;L^q) \quad \text{(for some $q>d$)},
\end{aligned}
\end{equation*}
then they coincide.\\
Moreover, the possible solution depends continuously on the initial conditions in the following sense: if $(v_i,\bC_i)$ is solution corresponding to the initial conditions $(v_{i,\mathrm{init}}, \bC_{i,\mathrm{init}})$, $i=\{1,2\}$, then there exists a constant $C\geq 0$ such that, for all $t\in (0,T)$,
$$
\|v_2(t)-v_1(t)\|_{L^2} + \|\bC_2(t)-\bC_1(t)\|_{L^2} \leq C \Big(  \|v_{2,\mathrm{init}} - v_{1,\mathrm{init}}\|_{L^2} + \|\bC_{2,\mathrm{init}}-\bC_{1,\mathrm{init}}\|_{L^2} \Big).
$$
\end{proposition}
\begin{proposition}[local existence]\label{th:local-existence1}
Let $d\in \{2,3\}$, $r\in ]1,+\infty[$ and $q\in]d,+\infty[$.\par\noindent
We assume that~$\Lambda$ is a positive function of class~$\mathcal C^1$ on~$\mathbb R_+$.\par\noindent
If $v_{\text{init}} \in D^r_q$ and~$\bC_{\text{init}}\in W^{1,q}$ is symmetric positive definite almost everywhere then there exists $T_\star>0$ and a strong solution $(v,\bC)$ to the system~\eqref{MGI0} in $[0,T_\star]$ satisfying:
\begin{equation*}
\begin{array}{ll}
 v \in L^r(0,T_\star;W^{2,q}),
& \partial_t v \in L^r(0,T_\star;L^q),\\
 \bC \in L^\infty(0,T_\star;W^{1,q}),
& \partial_t \bC \in L^r(0,T_\star;L^q).
\end{array}
\end{equation*}
\end{proposition}
\begin{theorem}[global existence]\label{th:global-existence1}
Under the same assumptions that in Proposition~\ref{th:local-existence1}, if we assume that
\begin{itemize}
\item $\displaystyle \frac{d}{2q}+\frac{1}{r}<\frac{1}{2}$,
\item $\displaystyle d=2 \quad \text{OR} \quad \big( d=3 \text{ but without taking into account the convective term $v\cdot \nabla v$} \big)$,
\end{itemize}
then the solution given by Proposition~\ref{th:local-existence1} holds taking any $T^\star >0$.
\end{theorem}

\subsubsection{Proofs of the MGI model results}

Recall that the MGI model~\eqref{MGI0} is based on the following constitutive relation, expressing the stress~$\bsigma$ as follows:
\[\bsigma = \Lambda(t) \frac{\sqrt{\bC}}{\tr \sqrt{\bC}} \quad \text{where} \quad \overset{\bigtriangledown}{\bC} = 0.\]
Since this model makes the root of~$\bC$ appear, we first note that the tensor~$\bC$ remains symmetric positive definite if it is initially the case (apply  Lemma~\ref{lem:conformation}). This observation suggests that the law is consistent.
In a first lemma (Lemma~\ref{lemma8}), we will see that it is possible to reformulate this constitutive relation into a differential equation linked~$\bsigma$, $\bsigma^2$ and~$Dv$.
Next, we will see in Lemma~\ref{lemma9} how to estimate $\nabla \bsigma$ using $\nabla (\bsigma^2)$.

\paragraph{Reformulation in term of stress tensor}

For the model~\eqref{MGI0}, it is possible to derive a differential equation for the stress tensor~$\bsigma$. The next idea follows the same kind of formulation that those presented in~\cite{Marrucci01}.

\begin{lemma}\label{lemma8}
Let $\Lambda:(0,T)\to \mathbb R$ be a positive function of class~$\mathcal C^1$.\par
If $\bsigma$ is defined by
\[\bsigma = \Lambda(t) \frac{\sqrt{\bC}}{\tr \sqrt{\bC}} \quad \text{where} \quad \overset{\bigtriangledown}{\bC} = 0,\]
then we have
\begin{equation}\label{MGI6}
\begin{aligned}
\overset{\bigtriangledown}{(\bsigma^2)} - \frac{2\Lambda'}{\Lambda}\bsigma^2 + \frac{2}{\Lambda} (Dv:\bsigma) \bsigma^2 = 0.
\end{aligned}
\end{equation}
\end{lemma}

\begin{remark}
In the two dimensional case, the Cayley-Hamilton theorem implies that
\[\bsigma^2 - (\tr \bsigma) \bsigma + (\det \bsigma) \bdelta = 0.\]
Since $\tr \bsigma =\Lambda(t)$ and $Dv:\bdelta = \div v = 0$, we have $\Lambda(t) Dv:\bsigma = Dv:\bsigma^2$. Consequently, the equation~\eqref{MGI6} can be written in terms of $\bsigma^2$ without intervening~$\bsigma$. From a numerical point of view, this remark can make it possible to treat this model almost as simply as the PEC model previously studied.
\end{remark}

{\bf Proof of Lemma~\ref{lemma8} -}
We introduce $\bB = \sqrt{\bC}$ and $\bQ = \frac{\bB}{\tr \bB}$. Derivating $\bQ^2$ with respect to the time (using precisely the convective derivative $d_t$), we obtain
\begin{equation}\label{MGI7}
\begin{aligned}
d_t (\bQ^2) 
 = d_t \Big( \frac{\bB^2}{(\tr \bB)^2} \Big) = \frac{d_t (\bB^2)}{(\tr \bB)^2} - 2 \frac{\tr (d_t \bB)}{(\tr \bB)^3} \bB^2.
\end{aligned}
\end{equation}
Writing the relation $\overset{\bigtriangledown}{\bC} = 0$ using the tensor $\bB$, we have
\begin{equation}\label{MGI8}
d_t (\bB^2) = \bB^2\cdot \nabla v + \transp{(\nabla v)} \cdot \bB^2.
\end{equation}
Multiplying left this equation by $\bB^{-1}$ and taking the trace, we obtain
\begin{equation}\label{MGI9}
d_t (\tr \bB) = Dv:\bB.
\end{equation}
Finally, plugging~\eqref{MGI8} and~\eqref{MGI9} in~\eqref{MGI7}, we obtain
\begin{equation*}
\begin{aligned}
d_t (\bQ^2) 
 = \frac{\bB^2\cdot \nabla v + \transp{(\nabla v)} \cdot \bB^2}{(\tr \bB)^2} - 2 \frac{Dv:\bB}{(\tr \bB)^3} \bB^2,
\end{aligned}
\end{equation*}
which can be completely rewritten in term of stress tensor~$\bQ$ as follows
\begin{equation*}
\begin{aligned}
\overset{\bigtriangledown}{(\bQ^2)} + 2 (Dv:\bQ) \bQ^2 = 0.
\end{aligned}
\end{equation*}
Finally, multiplying this equation by $\Lambda(t)^2$, we deduce the result since $\bsigma^2 = \Lambda(t)^2 \bQ^2$ with $\Lambda(t)>0$.
\finpreuve

\paragraph{Gradient of the square root}

\begin{lemma}\label{lemma9}
Let $\bsigma:\T^d \longrightarrow \mathbb R^{d\times d}$ be a symmetric positive definite matrix values function.\par
If $\bsigma$ is of class $\mathcal C^1$ then we have the following inequality
\begin{equation}\label{MGI10}
|\nabla \bsigma| \leq \frac{1}{\sqrt 2} |\bsigma^{-1}| |\nabla (\bsigma^2)|.
\end{equation}
\end{lemma}

\begin{remark}~\\[-0.5cm]
\begin{enumerate}
\item If $f:\T^d \longmapsto \mathbb R$ is a positive scalar function, the result is obvious since $\nabla (f^2) = 2 f \nabla f$. In particular, inequality~\eqref{MGI10} is not optimal.
\item If $\bsigma:\T^d \longrightarrow \mathbb R^{d\times d}$ is not symmetric positive definite then the result is false.
Consider for instance $\bsigma(x) = \begin{pmatrix} 1 & \cos(2\pi x_1) \\ 0 & -1 \end{pmatrix}$. We have  $\bsigma^2=\bdelta$ so $\nabla (\bsigma^2)=\bzero$ whereas $\nabla \bsigma \neq \bzero$.
\end{enumerate}
\end{remark}

{\bf Proof of Lemma~\ref{lemma9} - }
The spatial gradient of a $2$-tensor tensor is a $3$-tensor. We have
\[(\nabla (\bsigma^2))_{ijk} = \partial_i(\bsigma_{jm}\bsigma_{mk}) = \partial_i \bsigma_{jm} \bsigma_{mk} + \bsigma_{jm} \partial_i \bsigma_{mk}.\]
Using the notations given in the Appendix, page~\pageref{sec:notations}, it follows that
\[\transp{(\nabla (\bsigma^2))} = \transp{(\nabla \bsigma)} \cdot \bsigma + \bsigma \cdot \transp{(\nabla \bsigma)}.\]
Multiplying left and right by $\bsigma^{-1/2}$ (precisely making a $1$-contraction on left and right) we have
\begin{equation}\label{MGI11}
\bsigma^{-1/2} \cdot \transp{(\nabla (\bsigma^2))} \cdot \bsigma^{-1/2} = \underbrace{\bsigma^{-1/2} \cdot \transp{(\nabla \bsigma)} \cdot \bsigma^{1/2}}_{=\mathcal A} ~ + ~ \underbrace{\bsigma^{1/2} \cdot \transp{(\nabla \bsigma)} \cdot \bsigma^{-1/2}}_{=\mathcal B}.
\end{equation}
Computing the scalar product of the two $3$-tensors~$\mathcal A$ and~$\mathcal B$, we remark that
\[
\begin{aligned}
\mathcal A \overset{(3)}{:} \mathcal B 
& = \bsigma^{-1/2}_{im} ~ \transp{(\nabla \bsigma)}_{mjn} ~ \bsigma^{1/2}_{nk} ~~ \bsigma^{1/2}_{i\ell} ~ \transp{(\nabla \bsigma)}_{\ell jp} ~ \bsigma^{-1/2}_{pk} \\
& = \underbrace{\bsigma^{-1/2}_{im} \bsigma^{1/2}_{i\ell}}_{= \delta_{\ell m}} ~ (\nabla \bsigma)_{jmn} ~ (\nabla \bsigma)_{j\ell p} ~ \underbrace{\bsigma^{1/2}_{nk} \bsigma^{-1/2}_{pk}}_{= \delta_{np}},
\end{aligned}
\]
where we use the fact that $\bsigma$ is symmetric. We deduce that
$\mathcal A \overset{(3)}{:} \mathcal B = |\nabla \bsigma|^2$.
The equality~\eqref{MGI11} is written
\[| \bsigma^{-1/2} \cdot \transp{(\nabla (\bsigma^2))} \cdot \bsigma^{-1/2} |^2 = |\mathcal A|^2 + 2 \mathcal A \overset{(3)}{:} \mathcal B + |\mathcal B|^2 \geq 2 |\nabla \bsigma|^2,\]
which directly implies the desired result.
\finpreuve

\paragraph{Uniqueness and stability --}
The procedure to prove uniqueness and stability results is similar to the case of the PEC model (see the subsection~\ref{subsec:PEC-local-existence}): we write the difference between two solutions indexed by~$1$ and by~$2$, and then we make a $L^2$ estimate making the scalar product of the results by $v=v_1-v_2$ and~$\bC=\bC_1-\bC_2$.
The only difference with the case presented in the subsection~\ref{subsec:PEC-local-existence} comes from how to control the difference between the stresses $\bsigma=\bsigma_1-\bsigma_2$ from the difference in the conformations~$\bC$.

We recall that
\[\bsigma = \Lambda(t) \Big( \frac{\sqrt{\bC_1}}{\tr \sqrt{\bC_1}} - \frac{\sqrt{\bC_2}}{\tr \sqrt{\bC_2}} \Big).\]
The desired control is a consequence of the mean value theorem apply to the function
\[f:\bC\in \mathcal S_{+} \longmapsto \frac{\sqrt{\bC}}{\tr \sqrt{\bC}} \in \mathcal S_{+},\]
the set~$\mathcal S_{+}$ denoting the set of symmetric positive definite matrices.
Indeed, using the inverse function theorem, it's well known that the inverse of $\bC\in \mathcal S_{+} \longmapsto \bC^2 \in \mathcal S_{+}$ (that is the square root matrix function) is of class~$\mathcal C^1$. The function~$f$ is also of class $\mathcal C^1$ and the mean value theorem reads
\[|\bsigma| = |\Lambda(t)|\, |f(\bC_1) - f(\bC_2)| \leq \sup_{(0,T)} |\Lambda| \, M\, |\bC_1-\bC_2| \lesssim |\bC|,\]
where~$M$ is a bound on~$\mathrm df$ on any compact set of~$\mathcal S_{+}$ containing~$\bC_1$ and~$\bC_2$.

\paragraph{Local existence --}
The proof follows the same steps that in the proof of the local existence for the PEC type model (see the subsection~\ref{subsec:PEC-local-existence}).
We rewrite the system~\eqref{MGI0} as a fixed point equation
\begin{equation*}\label{MGI20}
\left\{
\begin{aligned}
& \partial_t v + \overline v\cdot\nabla \overline v - \Delta v + \nabla p = \div \bsigma \\
& \div v = 0 \\
& \bsigma =  \Lambda(t) \frac{\sqrt{\overline \bC}}{\tr \sqrt{\overline \bC}} \\
& \partial_t \bC + \overline v \cdot \nabla \bC - \bC\cdot \nabla \overline v - \transp{(\nabla \overline v)}\cdot \bC  = \bzero,
\end{aligned}
\right.
\end{equation*}
and apply Schauder’s theorem to the application $(\overline v, \overline \bC) \longmapsto (v,\bC)$ exactly like in the subsection~\ref{subsec:PEC-local-existence}.

\paragraph{Global existence --}

The strategy to show the result of global existence in time - Theorem~\ref{th:global-existence1} - follows the same principles as that of Theorem~\ref{th:global-existence}. The $L^\infty$-estimate on the stress~$\bsigma$ is a consequence of the estimate on the trace~$\tr \bsigma$ (the trace controls the $L^\infty$-norm for any symmetric positive definite matrix, see the proof of Lemma~\ref{lem:stress})~: \[|\bsigma(t,x)| \leq \tr \bsigma(t,x) = \Lambda(t) \leq \sup_{(0,T)}|\Lambda|.\] Lemma~\ref{lem:velocity} is again used to control the velocity from the stress.

To conclude the proof, it is enough to know how to control $y(t) = \|\nabla \bsigma\|_{L^r(0,t;L^q)}$ for $0\leq t\leq T_\star$.
In practice, we have no equation on $\bsigma$ but only one coupling $\bsigma^2$ and $\bsigma$ (see the equation~\eqref{MGI6} in Lemma~\ref{lemma8}). We then first get an estimate on $w(t) = \|\nabla (\bsigma^2)\|_{L^r(0,t;L^q)}$. We have
\begin{equation}\label{MGI21}
\begin{aligned}
\overset{\bigtriangledown}{(\bsigma^2)} -\frac{2 \Lambda'}{\Lambda} \bsigma^2 + \frac{2}{\Lambda} (Dv:\bsigma) \bsigma^2 = 0.
\end{aligned}
\end{equation}

We first take gradients in~\eqref{MGI21} and compute the $L^2$ scalar product of the resulting equation with~$q|\nabla (\bsigma^2)|^{q-2} \nabla (\bsigma^2)$. Using the Hölder inequality and the fact that $|\bsigma(t,\cdot)|_{L^\infty} \lesssim 1$ we deduce 
\begin{equation}\label{MGI22}
d_t \|\nabla (\bsigma^2)\|_{L^q}^q \lesssim (1 + \|\nabla v\|_{L^\infty}) \|\nabla (\bsigma^2)\|_{L^q}^q + \|\nabla^2 v\|_{L^q} \|\nabla (\bsigma^2)\|_{L^q}^{q-1} + \|\nabla v\|_{L^\infty} \|\nabla \bsigma\|_{L^q} \|\nabla (\bsigma^2)\|_{L^q}^{q-1}.
\end{equation}
Since $\bsigma$ is positive definite on a compact set $[0,T_\star]\times \T^d$, Lemma~\ref{lemma9} implies
\begin{equation}\label{MGI23}
\|\nabla \bsigma\|_{L^q} \lesssim \|\nabla (\bsigma^2)\|_{L^q}.
\end{equation}
The estimate~\eqref{MGI22} becomes
\begin{equation}\label{MGI24}
d_t \|\nabla (\bsigma^2)\|_{L^q}^q \lesssim (1 + \|\nabla v\|_{L^\infty}) \|\nabla (\bsigma^2)\|_{L^q}^q + \|\nabla^2 v\|_{L^q} \|\nabla (\bsigma^2)\|_{L^q}^{q-1}.
\end{equation}
The end of the proof is similar to the proof of the PEC case: the function $w(t) = \|\nabla (\bsigma^2)\|_{L^r(0,t;L^q)}$ satisfies
$$w'(t) \lesssim w(t) \ln (\mathrm e + w(t)) + w(t).$$
It is bounded by a double exponential function of kind $\mathrm e^{\mathrm e^t}$.
In virtue of~\eqref{MGI23}, that we can also write $y\lesssim w$, we deduce that~$y$ is bounded too, that concludes the proof.

\subsection{The Pom-Pom polymer model with constant backbone stretch}\label{subsec:pom-pom}

The Pom-Pom model, introduced by McLeish and Larson~\cite{Larson98} in 1998, was a breakthrough in the field of viscoelastic constitutive equations.
The model is developed, mainly, for long-chain branched polymers. The multiple branched molecule can be broken down into several individual modes. Each  mode  is  represented  by  a  backbone  between  two branch points, with a number of dangling arms on every end. The backbone is confined by a tube formed by other backbones, for details refer to~\cite{Larson98}. In the case where the backbone stretch, denoted by~$\Lambda$, is assumed to be constant (or is assumed to be known as a given function depending on time,~$\Lambda(t)$), the constitutive equation reads
\begin{equation*}
 \bsigma = \Lambda(t) \frac{\bC}{\tr \bC} \quad \text{where} \quad \overset{\bigtriangledown}{\bC} + (\bC - \frac{1}{d} \bdelta) = 0.
\end{equation*}
This model resembles the MGI model studied in Subsection~\ref{subsec:MGI} with two differences:
\begin{enumerate}
\item The present model does not appear the square root of~$\bC$;
\item The evolution of the tensor~$\bC$ seems a little bit more complicated...
\end{enumerate}
The first point simplifies the study since we do not need to control~$\sqrt{\bC}$ from~$\bC$ (that correspond to Lemma~\ref{lemma9}).
We also make sure that the second point does not disturb the proofs of Subsection~\ref{subsec:MGI}. In particular, we use the first part of the Lemma~\ref{lem:conformation} in order to verify that~$\bC$ still remains symmetric positive definite.
\begin{remark}
The original differential form by McLeish and Larson (as for the complete MGI model) make appear a differential equation implementing the evolution of the backbone stretch~$\Lambda$:
\begin{equation*}
\left\{
\begin{aligned}
& \bsigma = \Lambda^2 \frac{\bC}{\tr \bC} \\
& \overset{\bigtriangledown}{\bC} + (\bC - \frac{1}{d} \bdelta) = 0 \\
& d_t \Lambda = \Lambda (Dv:\bC) - (\Lambda-1) \quad \text{strictly for $\Lambda < q$}.
\end{aligned}
\right.
\end{equation*}
But the presence of the equation on~$\Lambda$ make appear mathematical difficulties: it seems difficult to have a $L^\infty$-estimates on~$\Lambda$, despite condition~$\Lambda < q$, which does not seem entirely compatible with~$\Lambda$'s evolution equation. Moreover, this lack of control is noticeable in numerical simulations. For instance, on page 280 of~\cite{Keunings01}, it is written:\\[-0.1cm]

\hspace{1cm}\parbox{0.85\linewidth}{
{\it "In the complex flow simulations, we also observed an unphysical behaviour for the pom–pom stretch parameter that cannot be observed in steady and start-up rheometrical flows. Indeed, at high Weissenberg number, the stretch may become smaller than unity when $Dv : \bC$ becomes negative, which is physically unrealistic. This happens in flow regions where the velocity gradient changes sign as, for example down- stream of a contraction/expansion."}
}\\

To avoid regularity issues when~$\Lambda$ reaches $q$, Renardy~\cite[page 97]{Renardy09} regularize the equation for~$\Lambda$~:
\[d_t \Lambda = \frac{q-\Lambda}{q-\Lambda+\eps} \Big[\Lambda (Dv:\bC) - (\Lambda-1) \Big],\]
with $\eps>0$ so that~$\Lambda<q$, and using similar method that for the MGI model, we can prove that~$\bsigma$ remains bounded. 
Since then, several authors have proposed relatively complex modifications, see for instance~\cite{Verbeeten01}. It would certainly be interesting to try to analyze these models from a mathematical point of view.
\end{remark}

\subsection{Example of "bad" model}\label{sec:false}

Several models take similar forms to those studied in this article.
This is for example the case of the non stretching form of the Rolie-Poly model which is written (see for instance~\cite{Renardy-Renardy17b}):
\begin{equation}\label{eq:09O1}
\left\{
\begin{aligned}
& \bsigma = G \bC \\
& \overset{\bigtriangledown}{\bC} + \frac{2}{3} (Dv:\bC)(\bC+\mu (\bC-\bdelta)) + \eps(\bC-\bdelta) = \bzero.
\end{aligned}
\right.
\end{equation}
\begin{remark}
We can notice that the solution of the model~\eqref{eq:09O1} satisfies $\tr(\bC-\bdelta)=0$ if this relationship is true initially.
However, this is not fine enough to apply the results proved in the present article since the trace control does not allow to have a bound on the tensor (especially if it is not defined positive!).
\end{remark}
Unfortunately, and unlike the other models presented in the present paper, it is not clear that the tensor~$\bC$, satisfying~\eqref{eq:09O1} always remains positive.
More precisely, the proof presented by Hulsen in~\cite{Hulsen90} gives a sufficient condition for the solution~$\bC$ of 
\begin{equation}\label{eq:cex1}
\overset{\bigtriangledown}{\bC} = g_1(\bC, \nabla v) \, \bdelta + g_2(\bC,\nabla v) \, \bC + g_3(\bC,\nabla v) \, \bC^2,
\end{equation}
to remain positive definite if the initial value $\bC_{\mathrm{init}}$ is positive definite.
In particular, the proof strongly use the assumption~$g_1>0$.
We can see that the result holds if $g_1 \equiv 0$ but it is possible to build a counterexample in the case where we have no information on the sign of~$g_1$, that is the case for model~\eqref{eq:09O1}: the value of $g_1$ being given by $g_1(\bC, \nabla v) = \varepsilon + \frac{2}{3}\mu Dv:\bC$.
\begin{example}
Consider the velocity field $v(x,y) = (x,-y)$ and the solution $\bC$ to the following equation
\begin{equation}\label{eq:cex2}
\overset{\bigtriangledown}{\bC} = (Dv:\bC)(\bdelta - \bC),
\end{equation}
corresponding to the equation~\eqref{eq:cex1} with $g_1(\bC, \nabla v) = Dv:\bC$, $g_2(\bC, \nabla v) = -Dv:\bC$ and $g_3=0$.\par
We assume that the initial condition is given by
\[\bC_{\mathrm{init}} = \begin{pmatrix} a_0  & 0 \\ 0 & b_0 \end{pmatrix}, \quad a_0>0, \quad b_0>0,\]
so that we can looking for a solution $\bC$ as a diagonal matrix:
\[\bC(t,x,y) = \begin{pmatrix} a(t)  & 0 \\ 0 & b(t) \end{pmatrix}.\]
The functions $a$ and $b$ satisfy
\begin{equation*}\label{eq:app2}
\left\{
\begin{aligned}
& a'-2a = 2(a-b)(1-a) \\
& b'+2b = 2(a-b)(1-b).
\end{aligned}
\right.
\end{equation*}
We look for $a_0>0$, $b_0>0$ such that $a(t)<0$ for~$t$ large enough.
Introducing $\alpha = a+b$ and $\beta=a-b$, we obtain the following equivalent system
\begin{equation*}\label{eq:app3}
\left\{
\begin{aligned}
& \alpha'-2\beta = 2\beta(2-\alpha) \\
& \beta'-2\alpha = -2 \beta^2.
\end{aligned}
\right.
\end{equation*}
In term of $(\alpha,\beta)$, we look for initial condition $(\alpha_0,\beta_0)$ with $-\alpha_0<\beta_0<\alpha_0$ and such that $\beta(t)<-\alpha(t)$ for~$t$ large enough.
Remark that if $\alpha_0=3$ then $\alpha(t)=3$, and in this case, $\beta$ satisfies
\begin{equation}\label{eq:app30}
\beta' + 2 \beta^2 = 6.
\end{equation}
We finally look for $\beta_0$ with $-3<\beta_0<3$ such that the solution of~\eqref{eq:app30} satisfies $\beta(t)<-3$ for~$t$ large enough.
Equation~\eqref{eq:app30} is a Riccati equation those the solution is expressed as
\[\beta(t) = \sqrt 3 - \frac{2\sqrt 3}{1+ \big( \frac{\sqrt 3 + \beta_0}{\sqrt 3-\beta_0} \big) \mathrm e^{4\sqrt 3 t}}.\]
If $\beta_0<-\sqrt 3$ then the function~$\beta$ decrease and goes to $-\infty$ in finite time: for~$t$ large enough, we have $\beta(t)<-3$.\\

As example, we take $\alpha_0=3$, $\beta_0=\frac{1}{2}(-3-\sqrt 3)$. In term of initial values $(a_0,b_0)$, we have $a_0 = \frac{1}{2}(\alpha_0+\beta_0) =\frac{3-\sqrt 3}{4}$ and $b_0 = \frac{1}{2}(\alpha_0-\beta_0) =\frac{9+\sqrt 3}{4}$.
We compute the time $t_1$ for which~$\beta(t_1)=-3$ (that corresponds to the time where $a(t_1) = \alpha(t_1) + \beta(t_1) = 0$):
\[t_1 = \frac{1}{8\sqrt 3}\ln 3.\]
As a consequence, if the initial condition is on the following form (which is a positive definite matrix)
\[\bC_{\mathrm{init}} = \frac{1}{4}\begin{pmatrix} 3-\sqrt 3  & 0 \\ 0 & 9+\sqrt 3 \end{pmatrix},\]
then the solution~$\bC$ of the equation~\eqref{eq:cex2} becomes non positive for time $t>\frac{1}{8\sqrt 3}\ln 3$.
\end{example}
From a physical point of view, to quote Leonov and Prokunin~\cite[Section 3.3.5, page 68]{Leonov94}:\\[-0.2cm]

\hspace{1cm}\parbox{0.85\linewidth}{
{\it "if in any point of the medium a principal value of~$\bC$ is negative, the Hadamard instability will immediately happen. Moreover, in this case, we cannot
also satisfy the restrictions imposed by the Second Law."}
}\\

In fact, such a behavior was numerically detected in~\cite[page 305]{Verbeeten02}:\\[-0.2cm]

\hspace{1cm}\parbox{0.85\linewidth}{
{\it "Physically, $\tr \bsigma$ can not become smaller than~$0$. However, numerically we
have encountered these unrealistic values at the front and back stagnation points in the flow around a cylinder at sufficiently high Weissenberg numbers."}
}

\section{Conclusion: How do you know if a rheological model is good?}\label{sec:conclusion}

In conclusion, the results proved in this paper show that adding non-linear contributions can generate fundamental and interesting mathematical properties. These types of results are well known when diffuse effects on stress are added, see for instance the works of Barrett and Suli~\cite{Barrett05, Barrett08, Barrett12, Barrett18}. Here we show that there are other ways to "regularize" that are hidden in relatively recent and relevant models: the partially extensible nature of the polymer strands for the PEC model, the chain stretch effects in convective constraint release theories for entangled polymers in the MGI model, etc.
The results of this article therefore affirm that adding such terms are not only essential to capture observed phenomena but also fundamental to show theoretical results: the two aspects naturally complement each other.\\
Another interesting point of view is that of the numerical aspect: numeric can be viewed as a bridge between the two aspects (experimental and theoretical) discussed in this conclusion. A mathematical model will be relevant if solutions can be calculated (at least approximately) and if these solutions are close to what is observed. But before calculating a solution numerically, it is essential to know if this solution exists, and if it is unique: that is if the model is mathematically well posed.
More precisely, in this article we add to the notion of well posed character in Hadamard's sense, the notion of the lifetime of the solution (because problems are evolutionary). We thus show that physically relevant problems for viscoelastic flow models are well posed in long time. These types of issues (globality in time of solutions) are generally not solved in similar frameworks, which generates numerical problems for evaluating solutions.
As example, the famous Height Weissenberg Number Problem (HWNP) is intimately linked to the lack of knowledge on the theoretical structure of the Oldroyd linear model.
Many authors published articles in the 1980s that explore the HWNP: Debbaut and Crochet in \cite{Debbaut}, Keunings in~\cite{Keunings0,Keunings1} or Davies in~\cite{Davies}. See also more recently the works~\cite{Boyaval09, Fattal04} without any real success. At the same time, many theoretical studies (see~\cite{Chemin,Fernandez-Guillen-Ortega, Guillope-Saut1, Lions-Masmoudi-viscoelastique}, but there are many others) are interested in the problems of global existence for the Oldroyd model, without more success!\\

Finally, the natural question can be: 
\begin{center}
Can we validate a model as soon as we know that it is mathematically well posed?
\end{center}
The results presented here tend to say that PEC and MGI models are mathematically relevant, and it would therefore be possible to calculate numerical solutions even in long time - with some validity.\\
The key to success in the models presented lies mainly on the natural bound that exists on stress. Roughly speaking, we show that if we have sufficient bounds on stress then the underlying models will be well posed. It is also important to note that for the study introduced here, these bounds are strongly related to a positiveness property (which is also physically important). The results of the present articles provide answers on some models but there are still many questions:\\[0.2cm]
\phantom{quad} $\bullet$ On the one hand, for some other models, the question of global existence remains open because we do not know how to obtain sufficient bounds on stress. This is the case with the Oldroyd model, although it is known that the associated stress is positive definite, see~\cite{LeBris12} (the positiveness is a consequence of the Kramers expression for the stress tensor using a micro-macro model). Beale–Kato–Majda criterion of explosion are known for viscoelastic models (see~\cite{Hu13,LeiMasmoudiZhou10})  but no result makes it possible to affirm if these criteria are or are not satisfied.\\[0.2cm]
\phantom{quad} $\bullet$ On the other hand, it seems clear that we can invalidate some models as soon as we have that they are not mathematically well posed. This is the subject of Subsection~\ref{sec:false}, which highlights some models for which the conformation tensor do not remain positive over time. The theoretical analysis proposed in this article can therefore no longer be applied, and the underlying physics therefore seems to be in serious trouble!

\section*{Appendix - Notations}\label{sec:notations}

In this appendix, we describe the notations used in the article. First, the notations on the tensor analysis are recalled (contracted products and norms). Then, we specify the notations used for the temporal derivatives (partial, total, convected...). Finally, the functional spaces used in the article are introduced.\\

First note that throughout the article, we use the symbol $x\lesssim y$ which means that there exists a constant~$C$ such that $x\leq Cy$. Obviously, according to the context, the constant~$C$ does not depend on the main variable of the studied problem. 

\paragraph{Tensorial tools}~\\[-0.3cm]

$\checkmark$ For a $2$-tensor~$\bA$, we denote by "$\tr \bA$" and by "$\det \bA$" its trace and its determinant respectively. The bold symbol~$\bdelta$ corresponds to the identity $2$-tensor whose components satisfies $\bdelta_{ij}=1$ if $i=j$ and $\bdelta_{ij}=0$ if not.\\[-0.2cm]

$\checkmark$ For a $N$-tensor, $N\geq 2$, we denote by $\transp{\bA}$ the generalization of the transpose operation. It is defined by inverting the two first indices:
\[\transp{\bA}_{ij{\bf k}} = \bA_{ji{\bf k}}.\]
In the previous expression - as in the following, the bold index notation~$\bf k$ denotes a multi-index. Above we have ${\bf k} = k_1,k_2,\cdots,k_{N-2}$.
As usual, we also introduce some contractions of pair of tensors. The $1$-contraction between a $N$-tensor~$\bA$ and a $M$-tensor~$\bB$ with $N\geq 1$, $M\geq 1$, is the $(N+M-2)$-tensor defined by (using the Einstein summation convention):
\[(\bA\cdot \bB)_{{\bf i}{\bf j}} = \bA_{{\bf i}m} \bB_{m{\bf j}}.\]
In the same way, the $2$-contraction between two tensors of order larger than~$2$ reads
\[(\bA : \bB)_{{\bf i}{\bf j}} = \bA_{{\bf i}mn} \bB_{mn{\bf j}}.\]
More generally, we can perform the $N$-contraction~$\overset{(N)}{:}$ of any pair of tensors of order larger than~$N$. For instance, we will use the $3$-contraction defined as follows
\[(\bA \overset{(3)}{:} \bB)_{{\bf i}{\bf j}} = \bA_{{\bf i}mn\ell} \bB_{mn\ell{\bf j}}.\]
$\checkmark$ Finally, the Frobenius norm of a $N$-tensor~$\bA$ will always denoted by $|\cdot|$, regardless of the size of~$N$. It is defined by
\[|\bA|^2 = \bA \overset{(N)}{:} \bA.\]

{\it
On several occasions, we will use the fact that when a $2$-tensor~$\bA$ is symmetrical and positive definite then its norm is equivalent to its trace:
\begin{equation}\label{eq:equivalent}
| \bA | \leq \tr \bA \leq \sqrt{d} \, |\bA|.
\end{equation}
The proof may be inferred from the orthogonal decomposition of symmetric positive definite $2$-tensor.
}

\paragraph{Time derivative}~\\[-0.3cm]

$\checkmark$ We denote the material derivative by~$d_t := \partial_t + v\cdot\nabla$.\\[-0.2cm]

$\checkmark$ It is well known that one must be careful when handling tensor equations because of the principle of frame indifference or material objectivity. The equations (and thus the time
derivative) should not depend on rigid motions of the observer.
We will use the classical upper convected time derivative (UCM) for a $2$-tensor function~$\bA$:
\begin{equation}\label{derive}
\overset{\bigtriangledown}{\bA} := d_t \bA - \bA\cdot \nabla v - \transp{(\nabla v)}\cdot \bA.
\end{equation}
{\it
Although we will only present results with the UCM derivative - except for an application to the PTT models in Subsection~\ref{subsec:PTT} , the theorems stated here can easily adapt to any classical frame indifferent derivatives~$\mathcal D_\xi$, $\xi\in [-1,1]$:
\begin{equation}\label{derive1}
\mathcal D_\xi \bA := d_t \bA + \frac{1-\xi}{2}( \bA\cdot \transp{\nabla v} + \nabla v \cdot \bA ) - \frac{1+\xi}{2} ( \bA \cdot \nabla v + \transp{(\nabla v)}\cdot \bA ),
\end{equation}
the particular cases $\xi=-1$, $\xi=0$ and $\xi=1$ correspond respectively to the lower convected, Jaumann and upper convected derivative - see Remarks~\ref{rem:1102}, point 3.
}

\paragraph{Functional spaces}~\\[-0.3cm]

$\checkmark$ The results presented in this article are only proved in space periodic domains~$\T^d$, $d=2$ or $d=3$. This constraint is purely technical and it is certainly possible that the results may extend to the case of the whole space~$\mathbb R^d$. Nevertheless, in the case of a bounded domain study, it must be ensured that all the arguments could still be adapted.\\[-0.2cm]

$\checkmark$ For all real $s\geq 0$ and all integer $q\geq 1$, the set $W^{s,q}$ corresponds to the Sobolev spaces. We classically denote $L^q=W^{0,q}$ the associated Lebesgue space.
Since we will frequently use functions with values in $\R^d$ or in the space~$\R^{d\times d}$ of real matrices... the usual notations will be abbreviated.
For instance, the space $(W^{1,q})^d$ will be denoted $W^{1,q}$.
Moreover, all functional norms will be denoted by indices, for instance like~$\| \cdot \|_{W^{1,q}}$.\\[-0.2cm]

$\checkmark$ The space $D^r_q$ stands for some fractional domain of the Stokes operator~$A_q$ in $L^q$ (cf. Section 2.3 in~\cite{Danchin}). Its norm is defined by
$$
\|w\|_{D^r_q} := \|w\|_{L^q} + \Big( \int_0^{+\infty} \|A_q \mathrm e^{-tA_q}w \|_{L^q}^{r} \, \mathrm dt \Big)^{1/r}.
$$
{\it
Roughly, the vector-fields of~$D^r_q$ are vectors which have $2-\frac{2}{r}$ derivatives in $L^q$ and are divergence-free.
It may be identified with Besov spaces. It also can be view as an interpolate space between $L^q$ and the domain of the Stokes operator~$D(A_q)$, see~\cite{Danchin}.
}\\[-0.2cm]

$\checkmark$ The notation of kind $L^r(0,T;W^{1,q})$ denotes the space of $r$-integrable functions on $(0,T)$, with values in the space~$W^{1,q}$.\\[-0.2cm]

$\checkmark$ Finally let us denote~$\mathds{P}$ the orthogonal projector in~$L^2$ onto the set of the divergence-free vectors fields of~$L^2$.
The pressure~$p$ is a Lagrange multiplier associated to the divergence free constraint.  
It can be solved using the Riesz transforms. More precisely, we take the divergence of the first equation of Navier-Stokes equations, and we use the periodic boundary conditions to have
\begin{equation}\label{pressure}
p=-(-\Delta)^{-1}\div\div\, (\bsigma-v\otimes v).
\end{equation}
We will see from Theorem~\ref{th:local-existence} or~\ref{th:local-existence1} that the solutions discussed in this paper satisfy $\bsigma - v \otimes v \in L^\infty(0,T;L^2)$. The pressure is meant to be given by~\eqref{pressure}.

\bibliographystyle{plain}
\bibliography{references}

\begin{thebibliography}{10}

\bibitem{Barrett05}
John~W. Barrett, Christoph Schwab, and Endre S\"uli.
\newblock Existence of global weak solutions for some polymeric flow models.
\newblock {\em Math. Models Methods Appl. Sci.}, 15(6):939--983, 2005.

\bibitem{Barrett08}
John~W. Barrett and Endre S\"uli.
\newblock Existence of global weak solutions to dumbbell models for dilute
  polymers with microscopic cut-off.
\newblock {\em Math. Models Methods Appl. Sci.}, 18(6):935--971, 2008.

\bibitem{Barrett12}
John~W. Barrett and Endre S\"uli.
\newblock Existence and equilibration of global weak solutions to kinetic
  models for dilute polymers {II}: {H}ookean-type models.
\newblock {\em Math. Models Methods Appl. Sci.}, 22(5):1150024, 84, 2012.

\bibitem{Barrett18}
John~W. Barrett and Endre S\"uli.
\newblock Existence of global weak solutions to the kinetic {H}ookean dumbbell
  model for incompressible dilute polymeric fluids.
\newblock {\em Nonlinear Anal. Real World Appl.}, 39:362--395, 2018.

\bibitem{KBKZ2}
B.~Bernstein, E.A. Kearsley, and L.J. Zapas.
\newblock A study of stress relaxation with finite strain.
\newblock {\em Trans. Soc. Rheol.}, 17:35--92, 1964.

\bibitem{Bird87}
R~Byron Bird, Robert~C Armstrong, and Ole Hassager.
\newblock Dynamics of polymeric liquids. vol. 1: Fluid mechanics.
\newblock 1987.

\bibitem{Bird77}
R.B. Bird, O.~Hassager, R.C. Armonstrong, and C.F. Curtiss.
\newblock {\em Dynamics of Polymeric Fluids}, volume~2 of {\em Kinetic Theory}.
\newblock John Wiley and Sons, New York, 1977.

\bibitem{Boyaval09}
S\'ebastien Boyaval, Tony Leli\`evre, and Claude Mangoubi.
\newblock Free-energy-dissipative schemes for the {O}ldroyd-{B} model.
\newblock {\em M2AN Math. Model. Numer. Anal.}, 43(3):523--561, 2009.

\bibitem{Brandon-Hrusa}
Deborah Brandon and William~J. Hrusa.
\newblock Global existence of smooth shearing motions of a nonlinear
  viscoelastic fluid.
\newblock {\em J. Integral Equations Appl.}, 2(3):333--351, 1990.

\bibitem{Chemin}
Jean-Yves Chemin and Nader Masmoudi.
\newblock About lifespan of regular solutions of equations related to
  viscoelastic fluids.
\newblock {\em SIAM J. Math. Anal.}, 33(1):84--112 (electronic), 2001.

\bibitem{Chilcott88}
MD~Chilcott and John~M Rallison.
\newblock Creeping flow of dilute polymer solutions past cylinders and spheres.
\newblock {\em Journal of Non-Newtonian Fluid Mechanics}, 29:381--432, 1988.

\bibitem{chupin13}
Laurent Chupin.
\newblock Existence results for the flow of viscoelastic fluids with an
  integral constitutive law.
\newblock {\em J. Math. Fluid Mech.}, 15(4):783--806, 2013.

\bibitem{Chupin14}
Laurent Chupin.
\newblock Global existence results for some viscoelastic models with an
  integral constitutive law.
\newblock {\em SIAM J. Math. Anal.}, 46(3):1859--1873, 2014.

\bibitem{Constantin07}
Peter Constantin.
\newblock Smoluchowski navier-stokes systems.
\newblock {\em Contemporary Mathematics}, 429:85, 2007.

\bibitem{Danchin}
R.~Danchin.
\newblock Density-dependent incompressible fluids in bounded domains.
\newblock {\em J. Math. Fluid Mech.}, 8(3):333--381, 2006.

\bibitem{Davies}
AR~Davies.
\newblock Reentrant corner singularities in non-newtonian flow. part i: Theory.
\newblock {\em Journal of Non-Newtonian Fluid Mechanics}, 29:269--293, 1988.

\bibitem{Debbaut}
B~Debbaut and MJ~Crochet.
\newblock Further results on the flow of a viscoelastic fluid through an abrupt
  contraction.
\newblock {\em Journal of Non-Newtonian Fluid Mechanics}, 20:173--185, 1986.

\bibitem{Doi88}
M.~Doi and S.F. Edwards.
\newblock {\em The theory of polymer dynamics}.
\newblock Oxford University Press, 1988.

\bibitem{Fattal04}
Raanan Fattal and Raz Kupferman.
\newblock Constitutive laws for the matrix-logarithm of the conformation
  tensor.
\newblock {\em Journal of Non-Newtonian Fluid Mechanics}, 123(2-3):281--285,
  2004.

\bibitem{Fernandez-Guillen-Ortega-CRAS}
E.~Fern{\'a}ndez-Cara, F.~Guill{\'e}n, and R.R. Ortega.
\newblock Existence et unicit\'e de solution forte locale en temps pour des
  fluides non newtoniens de type {O}ldroyd (version {$L\sp s$}--{$L\sp r$}).
\newblock {\em C. R. Acad. Sci. Paris S\'er. I Math.}, 319(4):411--416, 1994.

\bibitem{Fernandez-Guillen-Ortega}
E.~Fern{\'a}ndez-Cara, F.~Guill{\'e}n, and R.R. Ortega.
\newblock Some theoretical results concerning non-{N}ewtonian fluids of the
  {O}ldroyd kind.
\newblock {\em Ann. Scuola Norm. Sup. Pisa Cl. Sci. (4)}, 26(1):1--29, 1998.

\bibitem{Giga91}
Yoshikazu Giga and Hermann Sohr.
\newblock Abstract lp estimates for the cauchy problem with applications to the
  navier-stokes equations in exterior domains.
\newblock {\em Journal of functional analysis}, 102(1):72--94, 1991.

\bibitem{Guillope-Saut4}
C.~Guillop{\'e} and J.-C. Saut.
\newblock Global existence and one-dimensional nonlinear stability of shearing
  motions of viscoelastic fluids of {O}ldroyd type.
\newblock {\em RAIRO Mod\'el. Math. Anal. Num\'er.}, 24(3):369--401, 1990.

\bibitem{Guillope-Saut-CRAS}
C.~Guillop{\'e} and J.C. Saut.
\newblock R\'esultats d'existence pour des fluides visco\'elastiques \`a loi de
  comportement de type diff\'erentiel.
\newblock {\em C. R. Acad. Sci. Paris S\'er. I Math.}, 305(11):489--492, 1987.

\bibitem{Guillope-Saut3}
C.~Guillop{\'e} and J.C. Saut.
\newblock Existence results for the flow of viscoelastic fluids with a
  differential constitutive law.
\newblock {\em Nonlinear Anal.}, 15(9):849--869, 1990.

\bibitem{Guillope-Saut1}
C.~Guillop{\'e} and J.C. Saut.
\newblock Mathematical problems arising in differential models for viscoelastic
  fluids.
\newblock In {\em Mathematical topics in fluid mechanics (Lisbon, 1991)},
  volume 274 of {\em Pitman Res. Notes Math. Ser.}, pages 64--92. Longman Sci.
  Tech., Harlow, 1992.

\bibitem{Hadamard}
Jacques Hadamard.
\newblock Sur les probl\`emes aux dériv\'ees partielles et leur signification
  physique.
\newblock {\em Princeton University Bulletin}, 13:49--52, 1902.

\bibitem{Renardy2}
William~J. Hrusa and Michael Renardy.
\newblock A model equation for viscoelasticity with a strongly singular kernel.
\newblock {\em SIAM J. Math. Anal.}, 19(2):257--269, 1988.

\bibitem{Hu13}
Xianpeng Hu and Ryan Hynd.
\newblock A blowup criterion for ideal viscoelastic flow.
\newblock {\em J. Math. Fluid Mech.}, 15(3):431--437, 2013.

\bibitem{Hulsen90}
Martien~A Hulsen.
\newblock A sufficient condition for a positive definite configuration tensor
  in differential models.
\newblock {\em Journal of non-newtonian fluid mechanics}, 38(1):93--100, 1990.

\bibitem{Jourdain-Lebris-Lelievre}
Benjamin Jourdain, Tony Leli{\`e}vre, and Claude Le~Bris.
\newblock Existence of solution for a micro-macro model of polymeric fluid: the
  {FENE} model.
\newblock {\em J. Funct. Anal.}, 209(1):162--193, 2004.

\bibitem{Keunings0}
Roland Keunings.
\newblock On the high weissenberg number problem.
\newblock {\em Journal of Non-Newtonian Fluid Mechanics}, 20:209--226, 1986.

\bibitem{Keunings1}
Roland Keunings.
\newblock {\em Simulation of Viscoelastic Fluid Flow - Fundamentals of computer
  modeling for polymer processing}.
\newblock Hanser, 1989.

\bibitem{Kim}
Jong~Uhn Kim.
\newblock Global smooth solutions of the equations of motion of a nonlinear
  fluid with fading memory.
\newblock {\em Arch. Rational Mech. Anal.}, 79(2):97--130, 1982.

\bibitem{Larson84}
R.~G. Larson.
\newblock A constitutive equation for polymer melts based on partially
  extending strand convection.
\newblock {\em Journal of Rheology}, 28(5):545--571, 1984.

\bibitem{Larson97}
RG~Larson.
\newblock The elastic stress in “film fluids”.
\newblock {\em Journal of Rheology}, 41(2):365--372, 1997.

\bibitem{Larson13}
Ronald~G Larson.
\newblock {\em Constitutive Equations for Polymer Melts and Solutions:
  Butterworths Series in Chemical Engineering}.
\newblock Butterworth-Heinemann, 2013.

\bibitem{LeBris12}
Claude Le~Bris and Tony Leli{\`e}vre.
\newblock Micro-macro models for viscoelastic fluids: modelling, mathematics
  and numerics.
\newblock {\em Science China Mathematics}, 55(2):353--384, 2012.

\bibitem{LeiMasmoudiZhou10}
Zhen Lei, Nader Masmoudi, and Yi~Zhou.
\newblock Remarks on the blowup criteria for {O}ldroyd models.
\newblock {\em J. Differential Equations}, 248(2):328--341, 2010.

\bibitem{Lemarie}
P.~G. Lemari{\'e}-Rieusset.
\newblock {\em Recent developments in the {N}avier-{S}tokes problem}, volume
  431 of {\em Chapman \& Hall/CRC Research Notes in Mathematics}.
\newblock Chapman \& Hall/CRC, Boca Raton, FL, 2002.

\bibitem{Leonov94}
A.~I. Leonov and A.~N. Prokunin.
\newblock {\em Nonlinear Phenomena in Flows of Viscoelastic Polymer Fluids}.
\newblock Springer, Dordrecht, 1994.

\bibitem{Lions-Masmoudi-viscoelastique}
P.L. Lions and N.~Masmoudi.
\newblock Global solutions for some {O}ldroyd models of non-{N}ewtonian flows.
\newblock {\em Chinese Ann. Math. Ser. B}, 21(2):131--146, 2000.

\bibitem{Lodge64}
AS~Lodge.
\newblock Elastic liquids.
\newblock {\em London and New York}, pages 118--120, 1964.

\bibitem{Renardy17}
M\'aria Luk\'a\v~cov\'a Medvi\v~dov\'a, Hana Mizerov\'a, \v~S\'arka
  Ne\v~casov\'a, and Michael Renardy.
\newblock Global existence result for the generalized {P}eterlin viscoelastic
  model.
\newblock {\em SIAM J. Math. Anal.}, 49(4):2950--2964, 2017.

\bibitem{Marrucci00}
Giuseppe Marrucci, Francesco Greco, and Giovanni Ianniruberto.
\newblock Possible role of force balance on entanglements.
\newblock {\em Macromolecular Symposia}, 158(1):57--64, 8 2000.

\bibitem{Marrucci01}
Guiseppe Marrucci, Francesco Greco, and Giovanni Ianniruberto.
\newblock Integral and differential constitutive equations for entangled
  polymers with simple versions of ccr and force balance on entanglements.
\newblock {\em Rheologica Acta}, 40(2):98--103, Mar 2001.

\bibitem{Masmoudi08}
Nader Masmoudi.
\newblock Well-posedness for the {FENE} dumbbell model of polymeric flows.
\newblock {\em Comm. Pure Appl. Math.}, 61(12):1685--1714, 2008.

\bibitem{Masmoudi13}
Nader Masmoudi.
\newblock Global existence of weak solutions to the {FENE} dumbbell model of
  polymeric flows.
\newblock {\em Invent. Math.}, 191(2):427--500, 2013.

\bibitem{Larson98}
T.~C.~B. McLeish and R.~G. Larson.
\newblock Molecular constitutive equations for a class of branched polymers:
  The pom-pom polymer.
\newblock {\em Journal of Rheology}, 42(1):81--110, 1998.

\bibitem{Mitsoulis13}
E.~Mitsoulis.
\newblock 50 years of the k-bkz constitutive relation for polymers.
\newblock {\em ISRN Polymer Science}, 2013:22 pages, 2013.

\bibitem{NW02}
V~Ngamaramvaranggul and MF~Webster.
\newblock Simulation of pressure-tooling wire-coating flow with
  phan-thien/tanner models.
\newblock {\em International journal for numerical methods in fluids},
  38(7):677--710, 2002.

\bibitem{Oldroyd50}
JG~Oldroyd.
\newblock On the formulation of rheological equations of state.
\newblock In {\em Proceedings of the Royal Society of London A: Mathematical,
  Physical and Engineering Sciences}, volume 200, pages 523--541. The Royal
  Society, 1950.

\bibitem{Ottinger12}
Hans~C {\"O}ttinger.
\newblock {\em Stochastic processes in polymeric fluids: tools and examples for
  developing simulation algorithms}.
\newblock Springer Science \& Business Media, 2012.

\bibitem{PhanThien78}
Nhan Phan-Thien.
\newblock A nonlinear network viscoelastic model.
\newblock {\em Journal of Rheology}, 22(3):259--283, 1978.

\bibitem{Porod66}
G~Porod.
\newblock Fredrickson ag-principles and applications of rheology, 1966.

\bibitem{Renardy3}
M.~Renardy.
\newblock Existence of slow steady flows of viscoelastic fluids of integral
  type.
\newblock {\em Z. Angew. Math. Mech.}, 68(4):T40--T44, 1988.

\bibitem{Renardy}
M.~Renardy.
\newblock An existence theorem for model equations resulting from kinetic
  theories of polymer solutions.
\newblock {\em SIAM J. Math. Anal.}, 22(2):313--327, 1991.

\bibitem{Renardy09}
Michael Renardy.
\newblock Global existence of solutions for shear flow of certain viscoelastic
  fluids.
\newblock {\em J. Math. Fluid Mech.}, 11(1):91--99, 2009.

\bibitem{Renardy1}
Michael Renardy, William~J. Hrusa, and John~A. Nohel.
\newblock {\em Mathematical problems in viscoelasticity}, volume~35 of {\em
  Pitman Monographs and Surveys in Pure and Applied Mathematics}.
\newblock Longman Scientific \& Technical, Harlow, 1987.

\bibitem{Renardy-Renardy17b}
Michael Renardy and Yuriko Renardy.
\newblock A singular perturbation study of the rolie-poly model.
\newblock {\em Submit to J. Non-Newtonian Fluid Mech. - arXiv:1705.04677},
  2017.

\bibitem{Rivlin65}
RS~Rivlin.
\newblock Viscoelastic fluids.
\newblock {\em Research Frontiers in Fluid Dynamics}, page 144, 1965.

\bibitem{Saut13}
J.C. Saut.
\newblock Lectures on the mathematical theory of viscoelastic fluids.
\newblock {\em Lectures on the Analysis of Nonlinear Partial Differential
  Equations, Morningside Lectures in Mathematics}, (3):325--393, 2013.

\bibitem{Tanner00}
Roger~I Tanner.
\newblock {\em Engineering rheology}, volume~52.
\newblock OUP Oxford, 2000.

\bibitem{Tanner77}
Nhan~Phan Thien and Roger~I Tanner.
\newblock A new constitutive equation derived from network theory.
\newblock {\em Journal of Non-Newtonian Fluid Mechanics}, 2(4):353--365, 1977.

\bibitem{PTT77}
Nhan~Phan Thien and Roger~I Tanner.
\newblock A new constitutive equation derived from network theory.
\newblock {\em Journal of Non-Newtonian Fluid Mechanics}, 2(4):353--365, 1977.

\bibitem{Verbeeten01}
Wilco~MH Verbeeten, Gerrit~WM Peters, and Frank~PT Baaijens.
\newblock Differential constitutive equations for polymer melts: the extended
  pom--pom model.
\newblock {\em Journal of Rheology}, 45(4):823--843, 2001.

\bibitem{Verbeeten02}
Wilco~MH Verbeeten, Gerrit~WM Peters, and Frank~PT Baaijens.
\newblock Viscoelastic analysis of complex polymer melt flows using the
  extended pom--pom model.
\newblock {\em Journal of Non-Newtonian fluid mechanics}, 108(1):301--326,
  2002.

\bibitem{Keunings01}
P.~Wapperom and R.~Keunings.
\newblock Numerical simulation of branched polymer melts in transient complex
  flow using pom–pom models.
\newblock {\em Journal of Non-Newtonian Fluid Mechanics}, 97(2):267 -- 281,
  2001.

\bibitem{Zhang-Zhang}
Hui Zhang and Pingwen Zhang.
\newblock Local existence for the {FENE}-dumbbell model of polymeric fluids.
\newblock {\em Arch. Ration. Mech. Anal.}, 181(2):373--400, 2006.

\end{thebibliography}

\end{document}